\documentclass[final,supplement,onefignum,onetabnum]{siamart171218}

\usepackage{lipsum}
\usepackage{amsfonts}
\usepackage{graphicx}
\usepackage{epstopdf}
\usepackage{amsmath}

\usepackage{amsfonts}
\usepackage{amssymb}
\usepackage{algorithm}
\usepackage{algpseudocode}
\usepackage{listings}
\usepackage{booktabs}
\usepackage{epsfig}
\usepackage{color}
\usepackage{graphicx}
\usepackage{float}
\usepackage{hyperref}
\usepackage{footnote}
\usepackage{mathrsfs}
\usepackage{booktabs}
\usepackage{listings}
\usepackage{multirow}
\usepackage{caption}
\usepackage{subcaption}
\newcommand{\bpsi}{\boldsymbol\psi}
\newcommand{\B}{\mathcal{B}}

\newcommand{\divg}{\text{div}}
\newcommand{\relu}{\text{ReLU}}
\newcommand{\requ}{\text{ReQU}}

\newcommand{\normu}{\lVert u\rVert_{\B(\Omega)}}

\ifpdf
  \DeclareGraphicsExtensions{.eps,.pdf,.png,.jpg}
\else
  \DeclareGraphicsExtensions{.eps}
\fi

\newsiamremark{remark}{Remark}
\newsiamremark{hypothesis}{Hypothesis}
\crefname{hypothesis}{Hypothesis}{Hypotheses}
\newsiamthm{claim}{Claim}

\headers{Errors of deep mixed residual method}{L. Li, XC. Tai, J. Yang and Q. Zhu}

\title{Priori Error Estimate of Deep Mixed Residual Method for Elliptic PDEs\thanks{Submitted to the editors DATE.
\funding{This work is partially supported by the National Natural Science Foundation of China/Hong Kong RGC Joint Research Scheme (NSFC/RGC 11961160718). The work of J. Yang is supported by the National Science Foundation of China (NSFC-11871264).}}}

\author{Lingfeng Li\thanks{Department of Mathematics, Southern University of Science and Technology, Shenzhen, China;Department of Mathematics, Hong Kong Baptist University, Hong Kong, China (\email{lingfengli@life.hkbu.edu.hk})}
\and Xue-Cheng Tai\thanks{Department of Mathematics, Hong Kong Baptist University, Hong Kong, China (\email{xuechengtai@hkbu.edu.hk})}
\and Jiang Yang\thanks{Department of Mathematics, Southern University of Science and Technology, Shenzhen, China
(\email{yangj7@sustech.edu.cn}).}
\and Quanhui Zhu\thanks{Department of Mathematics, Southern University of Science and Technology, Shenzhen, China
(\email{12131244@mail.sustech.edu.cn}).}
}

\usepackage{amsopn}

\makeatletter
\newcommand*{\addFileDependency}[1]{% argument=file name and extension
  \typeout{(#1)}% latexmk will find this if $recorder=0 (however, in that case, it will ignore #1 if it is a .aux or .pdf file etc and it exists! if it doesn't exist, it will appear in the list of dependents regardless)
  \@addtofilelist{#1}% if you want it to appear in \listfiles, not really necessary and latexmk doesn't use this
  \IfFileExists{#1}{}{\typeout{No file #1.}}% latexmk will find this message if #1 doesn't exist (yet)
}
\makeatother

\newcommand*{\myexternaldocument}[1]{%
    \externaldocument{#1}%
    \addFileDependency{#1.tex}%
    \addFileDependency{#1.aux}%
}
\ifpdf
\hypersetup{
  pdftitle={Errors of deep mixed residual method},
  pdfauthor={L. Li, XC. Tai, J. Yang and Q. Zhu}
}
\fi

\myexternaldocument{ex_supplement}

\begin{document}

\maketitle

% REQUIRED
\begin{abstract}
In this work, we derive a priori error estimate of the mixed residual method when solving some elliptic PDEs. Our work is the first theoretical study of this method. We prove that the neural network solutions will converge if we increase the training samples and network size without any constraint on the ratio of training samples to the network size. Besides, our results suggest that the mixed residual method can recover high order derivatives better than the deep Ritz method, which has also been verified by our numerical experiments. 
\end{abstract}

% REQUIRED

% REQUIRED

\section{Introduction}
Using neural networks (NNs) to solve partial differential equations (PDEs) is very popular in these years. Since NNs are very powerful in approximating complicated functions, they are commonly used to solve PDEs in high dimension \cite{han2018solving,weinan2018deep,sirignano2018dgm,raissi2019physics,lu2021deepxde}. Besides, the neural networks can also be applied to approximate the solution operators \cite{lu2019deeponet,li2020fourier} and learn PDEs from data \cite{long2018pde,long2019pde}.

To approximate solutions to PDEs, we need first to define a subspace (subset) of functions and then search for a good solution in this subspace (subset). The classical finite element method (FEM) usually works with the finite elements space. To construct a finite element space, we need to discretize the domain and then define finite elements based on the mesh. The solution is then parameterized as a linear combination of all finite elements. FEM is very powerful in 2 and 3 dimensional problems. However, it suffers from the curse of dimensionality. Some mesh-free methods, e.g., the smoothed particle hydrodynamics \cite{gingold1977smoothed}, may be considered to solve high-dimensional problems, but these methods are very inaccurate near the boundary. NNs are used as a new mesh-free tool for PDEs in recent years. It is accurate for various problems, easy to implement and efficient to train. 

When solving PDEs with networks, we first parameterize the solution as a network with specified structures. Some commonly used networks are fully connected networks (FCN), residual networks (ResNet) and convolutional networks (CNN). A special case of the FCN is called two-layer network (TLN) which is often used in theoretical analysis \cite{barron1993universal,ma2019barron,lu2021priori}. The TLN is closely related to the Barron space \cite{barron1993universal} and we will also briefly introduce it in this paper. The training of PDE problems falls into the category of unsupervised training or semi-supervised training, where no ground truth or very few ground truth labels are known. 

Let us consider a general boundary value problem with the form:
\begin{equation*}
    \begin{aligned}
        f(u,\nabla u,\Delta u)&=0, &&x\in\Omega,\\
        B(u,\nabla u)&=0, &&x\in\partial\Omega,
    \end{aligned}
\end{equation*}
where $\Omega$ is a bounded domain. The DGM \cite{sirignano2018dgm} and PINN \cite{raissi2019physics} consider replacing $u$ with a neural network $\phi$ and then minimize the squared residual loss:
\begin{equation*}
    \Vert f(\phi,\nabla\phi,\Delta\phi)\Vert^2_{L^2(\Omega)}+\lambda\Vert B(\phi,\nabla\phi)\Vert^2_{L^2(\partial\Omega)},
\end{equation*}
where $\lambda$ is the penalty coefficient. Derivatives of $\phi$ with respect to $x$ are computed by the backpropagation. This type of method is simple and effective, but the evaluation of high order derivatives requires more computational efforts. To avoid the evaluation of high order derivatives, we may consider its variational form. If the equations is equivalent to a variational problem in a weak sense:
\begin{equation*}
    \min_{u} J(u):=a(u,u)-l(u),
\end{equation*}
where $a(\cdot,\cdot)$ is a bilinear form and $l$ is a linear functional, then we may use $J$ as the training loss and minimize $J(\phi)$. This is called deep Ritz method (DRM) \cite{weinan2018deep}. In some cases, we may add a penalty term for the boundary condition as well \cite{muller2021error}. \cite{zang2020weak} also proposed an adversarial network based on the variational form of PDEs. Another way to avoid higher order derivatives is to introduce an auxiliary vector-valued network $\psi$ to approximate $\nabla u$. Then, we combine the squared residual loss with an extra penalty term to form the mixed residual loss:
\begin{equation*}
    \Vert f(\phi,\bpsi,\divg(\bpsi))\Vert^2_{L^2(\Omega)}+\lambda_1\Vert B(\phi,\bpsi)\Vert^2_{L^2(\Omega)}+\lambda_2\Vert\nabla\phi-\bpsi\Vert^2_{L^2(\Omega)}.
\end{equation*}
This idea has been studied in \cite{lyu2022mim,cai2020deep,yang2021local}. \cite{lyu2020enforcing} further tried to enforce the boundary condition exactly to networks. This formulation is also called the first order system least square and has been applied to the finite element method \cite{cai1994first}. These loss functions represented in the form of integration are called expected loss functions. In practice, we would approximate these integrals with quadrature and this approximation is called the empirical loss function. 

In the study of numerical PDE methods, we are interested in the error estimates. However, for NNs, the error analysis is generally difficult, because of the non-linearity and non-convexity of the problem. The error of neural networks usually consists of two parts: approximation error and generalization error, which is also called quadrature error sometimes. The approximation error measures how well a given function can be approximated by a set of neural networks and the quadrature error measures the deviation between the expected loss and the empirical loss. So far, the error of DRM is well-studied in the literature \cite{lu2021priori,muller2021error,li2021generalization}. Like the finite element methods, the error $\Vert\phi-u^*\Vert_{H^1(\Omega)}$ is controlled by the $H^1$ approximation error of neural networks. \cite{mishra2020estimates,mishra2020estimates2} also gives a posteriori error for PINNs. For the Mixed residual method, there is no error analysis so far.

% \cite{shen2021deep} gives an explicit bound for the approximation error with respect to $L^p(\Omega)$ norm when $u^*\in C^s(\Omega)$. In the later work \cite{hon2021simultaneous}, the authors also generalize the approximation result to Sobolev norms. Though an explicit error bound is given, the error increase exponentially with the dimension $d$, which is undesired for high dimensional problems. In \cite{barron1993universal}, the authors derive a dimension-independent error estimates when approximating functions in the so-called Barron space by TLNs. The error is measured in $L^2(\Omega)$ norm and the rate only depends on the solution $u^*$ and the width of TLNs. Some later works \cite{ma2019barron,lu2021priori,li2021generalization} generalize the result to $H^1(\Omega)$ norm for functions in Barron-type spaces.

\textbf{Main contribution:} 
In this work, we provide a priori error analysis for the mixed residual method when solving some elliptic PDEs.
\begin{enumerate}
    \item This work is the first attempt at numerical analysis for the deep mixed residual method.
    \item Our result suggests that the mixed residual method can recover high order derivatives better than the deep Ritz method. We have also observed this property in numerical experiments.
    \item Our analysis covers both homogeneous and non-homogeneous boundary conditions while the analysis of deep Ritz method \cite{lu2021priori,muller2021error} only considers homogeneous cases.
    \item In this work, we estimate the approximation error and the quadrature error of two-layer networks equipped with a non-smooth activation function: rectified quadratic unit (ReQU) activation function.
    \item In \cite{lu2021priori}, the convergence of deep Ritz method requires that the number of training samples increases cubically with the network size. In our result, the increase of training samples is independent of network size.
\end{enumerate}
% In this work, we provide a priori error analysis for the mixed residual method when solving some elliptic PDEs.  We estimate the approximation error of sets of two-layer ReQU-activated networks using similar techniques with \cite{barron1993universal}. The quadrature error is estimated through the Rademacher complexity. When the true solution $u^*$ belongs to the Barron space, we can derive the convergence rate for $\Vert \phi-u^*\Vert_{H^1(\Omega)}$ and $\Vert \bpsi-\nabla u^*\Vert_{H_\divg}$ with respect to the networks size and number of quadrature points, i.e., training samples. Compared to the error estimates of DRM, our result suggests that the mixed residual method can approximate the Laplacian operator better, which has also been observed in numerical experiments as well. Besides, our error estimate works for non-homogeneous boundary conditions as well while the error estimates of DRM hold only for homogeneous boundary conditions. 

This paper is organized as follow. In Section \ref{sec:preliminaries}, we first introduce some preliminary results and describe our main theorems. In Section \ref{sec:error}, we will derive the error analysis for some elliptic PDEs. Then, we estimate the approximation error and quadrature error in Section \ref{sec:barron}. Lastly, we conduct some numerical experiments in Section \ref{sec:experiment} to partially verify our analysis.

\section{Preliminaries} \label{sec:preliminaries}
\subsection{Neural networks}
A neural network is a parametric function that maps an input $x\in\mathbb{R}^d$ to an output $y$. The output could be a scalar or a vector. Let us consider a simple FCN with scalar output:
\begin{equation*}\begin{aligned}
    h_j&=\sigma(W_{j-1}h_{j-1}+b_{j-1}),\quad j=2,\dots,L,\\
    y&=\omega^Th_L + b_L,
\end{aligned}\end{equation*}
where $h_0=x$, $\sigma$ is a non-linear activation function and $W_j$, $b_j$, $\omega$ are parameters to be learned. If the output $y$ is a vector, then we just replace $\omega$ with a matrix and replace $b_L$ with a vector. Another polular network structure is the ResNet \cite{he2016deep}:
\begin{equation*}\begin{aligned}
    h_j&=\sigma(W_{j-1}h_{j-1}+b_{j-1})+h_{j-1},\quad j=2,\dots,L,\\
    y&=\omega^Th_L + b_L.
\end{aligned}\end{equation*}
The main difference here is that an extra $h_{j-1}$ is added in each $h_j$ update, so the map $\sigma(W_{j-1}h_{j-1}+b_{j-1})$ equals to the residual $h_j-h_{j-1}$. This structure can avoid the so-called gradient vanishing problem in the training of very deep networks. One can refer to \cite{chen2018neural} for more theoretical explanation of ResNet.

Let us denote a network as $\phi(x)$. To train this network, we need to define an expected loss function $L(\phi)$ of the form:
\begin{equation*}
    L(\phi)=\int_{\Omega}l(x,\phi(x))\mathrm{d}\mu(x),
\end{equation*}
where $l(x,y)$ is a function that measures how well the network output $y=\phi(x)$ fits a given criterion, and $\mu(x)$ is a probability measure defined on $\Omega$. We then aim to minimize the loss $L(\phi)$ and find the optimal $\phi$ from a given set of networks: $\hat{\phi}=\arg\min_{\phi} L(\phi)$.
In practice, the integral form $L$ is usually approximated by a quadrature, i.e., we randomly sample some points from the distribution $\mu$ and approximate $L$ by
\begin{equation*}
    L_n(\phi)=\frac{1}{n}\sum_{i=1}^{n}l(x_i,\phi(x_i)).
\end{equation*}
$L_n$ is also called the empirical loss. The quadrature error is then defined as $|L_n(\phi)-L(\phi)|$.

\subsection{Deep mixed residual methods for elliptic PDEs}
Let $\Omega$ be a Lipschitz domain in the $d$-dimensional space. We consider the following second order PDEs with Neumann boundary condition
\begin{equation*}
    \begin{aligned}
        -\Delta u+ u=f \text{ for } x\in\Omega,\quad\text{and}\quad \frac{\partial u}{\partial \boldsymbol{n}}=g_1 \text{ for } x\in\partial\Omega, 
    \end{aligned}
    \label{eq:pde_neumann}
\end{equation*}
or Dirichlet boundary condition
\begin{equation*}
    \begin{aligned}
        -\Delta u=f \text{ for } x\in\Omega,\quad
        \text{and}\quad u=g_2 \text{ for } x\in\partial\Omega. 
    \end{aligned}
    \label{eq:pde_dirichlet}
\end{equation*}

To solve the PDE problems, we define two neural networks $\phi(x;\theta_1):\mathbb{R}^d\rightarrow\mathbb{R}$ and $\bpsi(x;\theta_2):\mathbb{R}^d\rightarrow\mathbb{R}^d$, where $\phi$ simulates the true solution $u^*$ and $\bpsi$ simulates $\nabla u^*$. Let $\phi\in V$ and $\bpsi\in W$, where $V$ and $W$ are sets of scalar-valued and vector-valued networks. The expected loss functions are given as below respectively. For the Neumann problem:
\begin{equation}
    L^N(\phi,\bpsi)=R_g(\phi,\bpsi)^2 + \lambda_1 R_{e}^N(\phi,\bpsi)^2 + \lambda_2 R_b^N(\phi,\bpsi),^2\label{eq:expected_loss_N}
\end{equation}
and for the Dirichlet problem:
\begin{equation}
    L^D(\phi,\bpsi)=R_g(\phi,\bpsi)^2 + \lambda_1 R_{e}^D(\phi,\bpsi)^2 + \lambda_2 R_b^D(\phi,\bpsi),^2\label{eq:expected_loss_D}
\end{equation}
where 
\begin{equation*}\begin{aligned}
    &R_g(\phi,\bpsi)=\Vert \nabla\phi-\bpsi\Vert_{L^2(\Omega)}, &\\
    &R_{e}^N(\phi,\bpsi)=\Vert -\divg(\bpsi)+\phi-f\Vert_{L^2(\Omega)},
    &&R_{e}^D(\phi,\bpsi)=\Vert -\divg(\bpsi)-f\Vert_{L^2(\Omega)},\\
    &R_b^N(\phi,\bpsi)=\Vert \bpsi\cdot\boldsymbol{n}-g_1\Vert_{L^2(\partial\Omega)},
    &&R_b^D(\phi,\bpsi)=\Vert \phi-g_2\Vert_{L^2(\partial\Omega)}.
\end{aligned}\end{equation*}
The corresponding empirical loss are respectively denoted by 
\begin{equation*}\begin{aligned}
    L_n^N(\phi,\bpsi)=&\frac{1}{n}\sum_{i=1}^n(|(\nabla\phi-\bpsi)(x_i)|^2+\lambda_1|(-\divg(\bpsi)+\phi-f)(x_i)|^2)\\
    &+\frac{1}{\bar{n}}\sum_{i=1}^{\bar{n}}\lambda_2|(\bpsi\cdot\mathbf{n}-g_1)(\bar{x}_i)|^2,
\end{aligned}\end{equation*}
and 
\begin{equation*}\begin{aligned}
    L_n^D(\phi,\bpsi)=&\frac{1}{n}\sum_{i=1}^n(|(\nabla\phi-\bpsi)(x_i)|^2+\lambda_1|(-\divg(\bpsi)-f)(x_i)|^2)\\&+\frac{1}{\bar{n}}\sum_{i=1}^{\bar{n}}\lambda_2|(\phi-g_2)(\bar{x}_i)|^2,
\end{aligned}\end{equation*}
where $\{x_i\}_{i=1}^n$ is randomly sampled from $\Omega$ and $\{\bar{x}_i\}_{i=1}^{\bar{n}}$ is the sample from $\partial\Omega$. In following discussions, $\bar{n}$ is an integer greater than or equals to $n/d^2$.
The error of neural network solutions are measured by $\Vert\phi-u^*\Vert_{H^1(\Omega)}^2$ and $\Vert\bpsi-\nabla u^*\Vert_{H_{\divg}(\Omega)}^2$ respectively.
% Here we briefly describe the training procedure of the deep mixed residual method in Algorithm \ref{alg:mixed}.

% \begin{algorithm}
% \caption{}\label{alg:mixed}
% \begin{algorithmic}[1]
% \Require A scalar-valued network $\phi:x\rightarrow\mathbb{R}$; A vector-valued network $\bpsi:x\rightarrow\mathbb{R}^d$; Total number of iterations $N_{itrs}$; The batch size $N_{bs}$; Two optimizers $O_1$ and $O_2$ for $\phi$ and $\bpsi$ respectively.
% \For{$i$ in range($N_{itrs}$)}
% \State Randomly sample $N_{bs}$ points in $\Omega$: $X_{int}$; $N_{bs}$ points in $\partial\Omega$: $X_{b}$
% \If{The boundary condition is Dirichlet condition}
% \State Evaluate $\phi(X_{int})$, $\phi(X_b)$, $\bpsi(X_{int})$
% \Else
% \State Evaluate $\phi(X_{int})$, $\bpsi(X_{int})$, $\bpsi(X_b)$
% \EndIf
% \State Compute the empirical loss function
% \State Update the optimizer $O_1$ and $O_2$ respectively

% \EndFor
% \end{algorithmic}
% \end{algorithm}

\subsection{Barron space and two-layer neural networks}
In this work, we assume that solutions to PDEs belong to the Barron space. Originally, the Barron space was introduced by \cite{barron1993universal}.
\begin{definition}
For a function $u$ defined on $\Omega$, if there exists a complex measure $F(d\omega)=e^{i\theta(\omega)}|\tilde{u}(\omega)|d\omega$ such that
\begin{equation*}
    u(x)=\int_{\mathbb{R}^d} e^{i\omega\cdot x}F(d\omega),\quad x\in\Omega,
\quad \text{and} \quad
    \int_{\mathbb{R}^d}(1+|\omega|_1)^s|\tilde{u}(\omega)|d\omega<+\infty,
\end{equation*}
then $u$ belongs to the Barron space $\B^s(\Omega)$ associated with the Barron norm 
\begin{equation*}
    \Vert u\Vert_{\B^s}:=\inf_{\tilde{u}}\int_{\mathbb{R}^d}(1+|\omega|_1)^s|\tilde{u}(\omega)|d\omega.
\end{equation*}
\end{definition}
The infimum in the definition is taken over all the possible complex measure $F(d\omega)=e^{i\theta(\omega)}|\tilde{u}(\omega)|d\omega$.
When $s=3$, we simply denote $\B^s$ as $\B$. It is easy to see that $\B^{s+1}\subset\B^{s}$ for $s\geq 1$ and \cite{siegel2020approximation} proves that $\B^s(\Omega)\subset H^s(\Omega)$ for any positive integer $s$ and bounded domain $\Omega\subset\mathbb{R}^d$.
% \begin{lemma}[\cite{siegel2020approximation}]
% $\B^s(\Omega)\subset H^s(\Omega)$ for any positive integer $s$ and bounded domain $\Omega\subset\mathbb{R}^d$.
% \end{lemma}

Here we also generalize the definition of Barron space to vector-valued functions.
\begin{definition}
For a vector-valued function $\boldsymbol{u}=(u_1,\dots,u_d)$ defined on $\Omega$, if $u_j\in\B^s(\Omega)$ for each $j=1,\dots,d$, we say $\boldsymbol{u}$ belongs to the Barron space $\B^s(\Omega;\mathbb{R}^d)$ associated with the Barron norm $\Vert \boldsymbol{u}\Vert_{\B^s(\Omega;\mathbb{R}^d)}=(\sum_{j=1}^d\Vert u_j\Vert_{\B^s(\Omega)}^2)^{1/2}$.
\end{definition}

\begin{lemma}
Given a function $u\in\B^{s+1}(\Omega)$, we have $\nabla u\in\B^s(\Omega;\mathbb{R}^d)$ and \\
$\Vert\nabla u\Vert_{\B^s(\Omega;\mathbb{R}^d)}\leq\Vert u\Vert_{\B^{s+1}(\Omega)}$.
\label{lem::vector_bnorm}
\end{lemma}
\begin{proof}
Given any $\eta>0$, there exist a $F(d\omega)=e^{i\theta(\omega)}|\tilde{u}(\omega)|d\omega$ such that $u(x)=\int_{\mathbb{R}^d}e^{i\omega\cdot x}F(d\omega)$ in $\Omega$ and $\int_{\mathbb{R}^d}(1+|\omega|_1)^{s+1}|\tilde{u}(\omega)|d\omega< \Vert u\Vert_{\B^{s+1}(\Omega)}+\eta$. Since the gradient of $u$ has the representation $\nabla u(x)=\int_{\mathbb{R}^d}e^{i\omega\cdot x}i\omega F(d\omega)$, then
\begin{equation*}\begin{aligned}
    \Vert\nabla u\Vert_{\B^s(\Omega;\mathbb{R}^d)}^2&\leq\bigg(\int_{\mathbb{R}^d}(1+|\omega|_1)^s|\omega|_1|\tilde{u}(\omega)|d\omega\bigg)^2\\
    &<\bigg(\int_{\mathbb{R}^d}(1+|\omega|_1)^{s+1}|\tilde{u}(\omega)|d\omega\bigg)^2<\bigg(\Vert u\Vert_{\B^{s+1}(\Omega)}+\eta\bigg)^2.
\end{aligned}\end{equation*}
Let $\eta$ go to $0$ and we complete the proof.
\end{proof}

For functions in the Barron space $\B(\Omega)$, we can use two-layer ReLU activated networks to approximate them to arbitrary actuary. \cite{li2021generalization} gives an estimation of the approximation error with the order $O(\frac{1}{\sqrt{m}})$, where $m$ is width of two-layer ReLU networks.
In our discussions, the activation function needs to be secondly differentiable. Instead, we use the rectified quadratic unit(ReQU) function $\requ(x) = \relu(x)^2=(\max\{x,0\})^2$.
We further define the set of ReQU activated networks associated with a Barron space function $u$ by
\begin{equation*}\begin{aligned}
    V^m_u=\bigg\{&c+\frac{1}{m}\sum_{i=1}^m a_i\requ(\omega_i^\top x+b_i)\bigg|\\
    &|c|\leq 2\Vert u\Vert_{\B}, |a_i|\leq 8\Vert u\Vert_{\B}, |\omega_i|_2\leq 1, |b_i|\leq 1, i=1,\dots,m \bigg\}.
    \label{eq::two_layer_nn}
\end{aligned}\end{equation*}
and
\begin{equation*}\begin{aligned}
    W^m_{\boldsymbol u}=\bigg\{&\boldsymbol{c}+\frac{1}{m}\sum_{i=1}^m \boldsymbol{a}_i\circ\requ(\boldsymbol{W}_i x+\boldsymbol{b}_i)\bigg|
    |\boldsymbol{c}_j|\leq 2\Vert \boldsymbol{u}_j\Vert_{\B}, 
    |(\boldsymbol{a}_i)_j|\leq 8\Vert \boldsymbol{u}_j\Vert_{\B}, \\
    &\|\boldsymbol{W}_i\|_2\leq 1, \|\boldsymbol{b}_i\|_\infty\leq 1, i=1,\dots,m, j=1,\dots,d \bigg\}.
\end{aligned}\end{equation*}

\subsection{Main results}
We first give our main results in this section. If we assume the domain $\Omega=[0,1]^d$, the errors of the mixed residual method for two PDE problems are given in the following theorems. Basically, the errors of neural networks can be decomposed into two parts: the approximation error and the quadrature error. 
\begin{theorem}\label{theorem:main_N}
Suppose the solution to Neumann problem (\ref{eq:pde_neumann}) $u^*_N$ belongs to the Barron space $\B^4(\Omega)$. Let $\hat{\phi},\hat{\bpsi}$ be the neural network solutions:
$$(\hat{\phi},\hat{\bpsi})=\mathop{\arg\min}\limits_{\phi\in V^m_{u^*_N},\bpsi\in W^m_{u^*_N}}L^N_n(\phi,\bpsi).$$
Then, we have
\begin{equation*}\begin{aligned}
    \mathbb{E}\Vert\hat{\phi}-u^*_N\Vert_{H^1(\Omega)}^2\lesssim \frac{\Vert u^*_N\Vert^2_{\B^4(\Omega)}}{m}+\frac{\Vert u^*_N\Vert^2_{\B^4(\Omega)}}{\sqrt{n}},
\end{aligned}\end{equation*}
and
\begin{equation*}\begin{aligned}
    \mathbb{E}\Vert\hat{\bpsi}-\nabla u^*_N\Vert_{H_{div}(\Omega)}^2\lesssim \frac{\Vert u^*_N\Vert^2_{\B^4(\Omega)}}{m}+\frac{\Vert u^*_N\Vert^2_{\B^4(\Omega)}}{\sqrt{n}},
\end{aligned}\end{equation*}
where the expectation is taken on the random sampling of training data in $\Omega$ and $\partial\Omega$.
\end{theorem}

\begin{theorem}\label{theorem:main_D}
Suppose the solution to Dirichlet problem (\ref{eq:pde_dirichlet}) $u^*_D$ belongs to the Barron space $\B^4(\Omega)$. Let $\hat{\phi},\hat{\bpsi}$ be the neural network solutions:
$$(\hat{\phi},\hat{\bpsi})=\mathop{\arg\min}\limits_{\phi\in V^m_{u^*_D},\bpsi\in W^m_{u^*_D}}L^D_n(\phi,\bpsi).$$
Then, we have
\begin{equation*}\begin{aligned}
    \mathbb{E}\Vert\hat{\phi}-u^*_D\Vert_{H^1(\Omega)}^2\lesssim \frac{\Vert u^*_D\Vert^2_{\B^4(\Omega)}}{\sqrt{m}}+\frac{\Vert u^*_D\Vert^2_{\B^4(\Omega)}}{n^{1/4}},
\end{aligned}\end{equation*}
and
\begin{equation*}\begin{aligned}
    \mathbb{E}\Vert\hat{\bpsi}-\nabla u^*_D\Vert_{H_{\divg}(\Omega)}^2\lesssim \frac{\Vert u^*_D\Vert^2_{\B^4(\Omega)}}{\sqrt{m}}+\frac{\Vert u^*_D\Vert^2_{\B^4(\Omega)}}{n^{1/4}},
\end{aligned}\end{equation*}
where the expectation is taken on the random sampling of training data.
\end{theorem}

\section{Error analysis for elliptic problems} \label{sec:error}
In this section, we will show how the error of neural networks can be controlled by the approximation error and quadrature error for elliptic PDEs with different boundary conditions. 
\subsection{Neumann boundary condition}
Suppose $V$ and $W$ be sets of scalar-valued and vector-valued neural networks with fixed structures and activation functions in $H^2(\mathbb{R})$. Then, for the PDE problem (\ref{eq:pde_neumann}) with Neuman boundary condition, we have the following error estimate.
\begin{theorem}\label{theorem:neumman}
Let $u^*_N$ be the classical solution to the Neumann problem (\ref{eq:pde_neumann}) and
% $$(\phi^*,\bpsi^*)=\underset{\phi\in V,\bpsi\in W}{\arg\min} L^N(\phi,\bpsi),$$
\begin{equation*}
    (\hat{\phi}, \hat{\bpsi})=\underset{\phi\in V,\bpsi\in W}{\arg\min} L_n^N(\phi,\bpsi).
\end{equation*}
Then, we have
\begin{equation*}
    \Vert\hat{\phi}-u^*_N\Vert^2_{H^1(\Omega)}\lesssim E_1+E_2,
    \quad \text{and} \quad
    \Vert\hat{\bpsi}-\nabla u^*_N\Vert_{H_{\divg}(\Omega)}^2\lesssim E_1+E_2,
\end{equation*}
where
$$E_1=2\sup_{\phi\in V,\bpsi\in W}|L^N(\phi,\bpsi)-L^N_n(\phi,\bpsi)|,$$
and
\begin{equation*}\begin{aligned}
E_2=\underset{\bpsi\in W}{\min}\left(\Vert\bpsi-\nabla u^*_N\Vert_{H_{\divg}(\Omega)}^2+\Vert \bpsi-\nabla u^*_N\Vert_{H^1(\Omega)}^2\right)+\underset{\phi\in V}{\min}\Vert \phi-u^*_N\Vert_{H^1(\Omega)}^2.
\end{aligned}\end{equation*}
\end{theorem}
% \begin{lemma}[Trace theorem]
% If the domain $\Omega$ has Lipschitz boundary, then there exist a continuous linear operator $\gamma:H^1(\Omega)\rightarrow L^2(\partial\Omega)$ that give functions in $C(\bar{\Omega})\cap H^1(\Omega)$ the classical boundary value, and there exist a constant $C_{\Omega}$ depends on the domain $\Omega$ such that
% \begin{equation*}
%     \Vert \gamma u\Vert_{L^2(\partial\Omega)}^2\leq C_{\Omega}\Vert u\Vert_{H^1(\Omega)}^2
% \end{equation*}
% for any $u\in H^1(\Omega)$.
% \end{lemma}
To prove Theorem \ref{theorem:neumman}, we first need to show the errors of two networks can be bounded by some multiples of the expected loss functions.
\begin{lemma}\label{lemma:posterror_1}
Let $\phi\in V$ and $\bpsi\in W$ be two networks, and $u^*_N$ be the classical solution to the Neumann problem (\ref{eq:expected_loss_N}). Then we have 
\begin{equation*}
\Vert\phi-u^*_N\Vert^2_{H^1(\Omega)}\leq 4C_{\Omega}(R^N_b)^2+2R_g^2+2(R^N_{e})^2,
\end{equation*}
and
\begin{equation*}
\Vert\bpsi-\nabla u^*_N\Vert_{H_{\divg}(\Omega)}^2\leq 6R_g^2+6(R_{e}^N)^2+8C_{\Omega}(R_b^N)^2.
\end{equation*}
\end{lemma}
\begin{proof}

Since $u^*_N$ is the classical solution to the PDE problem (\ref{eq:pde_neumann}), it satisfy the variational form:
\begin{equation*}
    \int_{\Omega}\nabla u^*_N\nabla v+u^*_Nv\ dx=\int_{\Omega}fv\ dx+\int_{\partial\Omega}g_1v\ ds,
\end{equation*}
for any $v\in H^1(\Omega)$. Let $v=\hat{\phi}=\phi-u^*_N$, we have
\begin{equation*}
    \int_{\Omega}\nabla u^*_N\nabla \hat{\phi}+u^*_N\hat{\phi}\ dx=\int_{\Omega}f\hat{\phi}\ dx+\int_{\partial\Omega}g_1\hat{\phi}\ ds.
\end{equation*}
Then, by adding $\int_{\Omega}(\divg(\bpsi)-\phi)\hat{\phi}\ dx-\int_{\partial\Omega}\hat{\phi}\nabla\phi\cdot\mathbf{n}\ ds$ on both sides,
\begin{equation*}\begin{aligned}
    &\int_{\Omega}(f+\divg(\bpsi)-\phi)\hat{\phi}\ dx-\int_{\partial\Omega}\hat{\phi}(\nabla\phi-\nabla u^*_N)\cdot \mathbf{n}\ ds\\
    =&\int_{\Omega}(\nabla u^*_N-\bpsi)\cdot\nabla \hat{\phi}+(u^*_N-\phi)\hat{\phi}\ dx-\int_{\partial\Omega}\hat{\phi}(\nabla\phi-\bpsi)\cdot \mathbf{n}\ ds\\
    =&\int_{\Omega}(\nabla u^*_N-\nabla\phi+\nabla\phi-\bpsi)\cdot\nabla \hat{\phi}+(u^*_N-\phi)\hat{\phi}\ dx+\int_{\partial\Omega}\hat{\phi}(\bpsi-\nabla\phi)\cdot \mathbf{n}\ ds\\
    =&-\int_{\Omega}|\nabla\hat{\phi}|^2\ dx-\int_{\Omega}|\hat{\phi}|^2\ dx+\int_{\Omega}(\nabla\phi-\bpsi)\cdot\nabla\hat{\phi}\ dx+\int_{\partial\Omega}\hat{\phi}(\bpsi-\nabla\phi)\cdot \mathbf{n}\ ds.
       \end{aligned}\end{equation*}
By rearranging the equation, we get
\begin{equation*}\begin{aligned}
    &\int_{\Omega}|\nabla\hat{\phi}|^2\ dx+\int_{\Omega}|\hat{\phi}|^2\ dx\\
    =&-\int_{\Omega}(f+\divg(\bpsi)-\phi)\hat{\phi}+\int_{\Omega}(\nabla\phi-\bpsi)\cdot\nabla\hat{\phi}\     dx+\int_{\partial\Omega}\hat{\phi}(\bpsi-\nabla u^*_N)\cdot \mathbf{n}\ ds\\
    =& \int_{\partial\Omega}\hat{\phi}(\bpsi\cdot\mathbf{n}-g_1)\ ds +\int_{\Omega}(\nabla\phi-\bpsi)\cdot\nabla\hat{\phi}\ dx+\int_{\Omega}(-\divg(\bpsi)+\phi-f)\hat{\phi}\ dx\\
   %\leq& C_{\Omega}\Vert\bpsi\cdot\mathbf{n}-g_1\Vert^2_{L^2(\partial\Omega)}+\frac{\Vert\hat{\phi}\Vert_{L^2(\partial\Omega)}^2}{4C_{\Omega}}+\frac{R_g^2}{2}+\frac{1}{2}\int_{\Omega}|\nabla\hat{\phi}|^2\ dx+\frac{(R^N_{e})^2}{2}+\frac{1}{2}\int_{\Omega}|\hat{\phi}|^2\ dx\\
    \leq& C_{\Omega}(R^N_b)^2+\frac{\Vert\hat{\phi}\Vert_{H^1(\Omega)}^2}{4}+\frac{R_g^2}{2}+\frac{1}{2}\int_{\Omega}|\nabla\hat{\phi}|^2\ dx+\frac{(R^N_{e})^2}{2}+\frac{1}{2}\int_{\Omega}|\hat{\phi}|^2\ dx,
\end{aligned}\end{equation*}
where the last inequality is derived by applying the trace theorem and $C_\Omega$ is a constant depends on the domain $\Omega$. Then, we combine like terms and get
\begin{equation*}
    \Vert\phi-u^*_N\Vert^2_{H^1(\Omega)}=\Vert\hat{\phi}\Vert_{H^1(\Omega)}^2\leq 4C_{\Omega}(R_b^N)^2+2R_g^2+2(R_{e}^N)^2.
\end{equation*}
For $\Vert\bpsi-\nabla u^*_N\Vert_{H_{\divg}(\Omega)}^2$, we have
\begin{equation*}\begin{aligned}
    &\Vert\bpsi-\nabla u^*_N\Vert_{H_{\divg}(\Omega)}^2\\
    %&=\Vert\bpsi-\nabla u^*_N\Vert_{L^2(\Omega)}^2+\Vert\divg(\bpsi)-\Delta u^*_N\Vert^2_{L^2(\Omega)}\\
    &=\Vert\bpsi-\nabla\phi+\nabla\phi-\nabla u^*_N\Vert_{L^2(\Omega)}^2+\Vert\divg(\bpsi)+f-\phi+\phi-u^*_N\Vert^2_{L^2(\Omega)}\\
    &\leq 2R_g^2+2(R_{e}^N)^2+2\Vert\hat{\phi}\Vert_{H^1(\Omega)}^2
    \leq 6R_g^2+6(R_{e}^N)^2+8C_{\Omega}(R_b^N)^2.
\end{aligned}\end{equation*}
\ 
\end{proof}

Next, we are going to prove the expected loss can be bounded by the approximation error of function classes $V$ and $W$ to the solution $u^*_N$. 
\begin{lemma}\label{lemma:losserror_1}
Let $ (\phi^*,\bpsi^*):=\underset{\phi\in V,\bpsi\in W}{\arg\min} L^N(\phi,\bpsi).$ Then,
\begin{equation*}\begin{aligned}
    L^N(\phi^*,\bpsi^*)\lesssim \underset{\bpsi\in W}{\min}\left(\Vert\bpsi-\nabla u^*_N\Vert_{H_{\divg}(\Omega)}^2+\Vert\bpsi-\nabla u^*_N\Vert_{H^1(\Omega)}^2\right)+\underset{\phi\in V}{\min}\Vert\phi-u^*_N\Vert_{H^1(\Omega)}^2.
\end{aligned}\end{equation*}
\end{lemma}
% \textbf{Remark. } In order to apply the trace theorem, we need to further assume that each component of $\nabla u^*_N$ belongs to $H^1(\Omega)$. If $u^*_N$ is assumed to be a Barron function, then this assumption is automatically satisfied.
\begin{proof}
For any $\tilde{\phi}\in V$ and $\tilde{\bpsi}\in W$, we have
\begin{equation*}\begin{aligned}
    L^N&(\phi^*,\bpsi^*)\leq L^N(\tilde{\phi},\tilde{\bpsi})\\
    =&\Vert \tilde{\bpsi}-\nabla\tilde{\phi}\Vert_{L^2(\Omega)}^2+\lambda_1\Vert-\divg(\tilde{\bpsi})+\tilde{\phi}-f\Vert_{L^2(\Omega)}^2+\lambda_2\Vert (\tilde{\bpsi}-\nabla u^*_N)\cdot\mathbf{n}\Vert_{L^2(\partial\Omega)}^2\\
    \leq&\Vert \tilde{\bpsi}-\nabla u^*_N+\nabla u^*_N-\nabla\tilde{\phi}\Vert_{L^2(\Omega)}^2+\lambda_2\Vert\tilde{\bpsi}-\nabla u^*_N\Vert_{L^2(\partial\Omega)}^2\\
    &+\lambda_1\Vert-\divg(\bpsi)+\divg(\nabla u^*_N)+\tilde{\phi}-u^*_N\Vert_{L^2(\Omega)}^2\\
    \leq& 2\Vert \tilde{\bpsi}-\nabla u^*_N\Vert_{L^2(\Omega)}^2+2\Vert \nabla\tilde{\phi}-\nabla u^*_N\Vert_{L^2(\Omega)}^2+\lambda_2C_{\Omega}\Vert\tilde{\bpsi}-\nabla u^*_N\Vert_{H^1(\Omega)}^2\\
    &+2\lambda_1\Vert\divg(\tilde{\bpsi})-\divg(\nabla u^*_N)\Vert_{L^2(\Omega)}^2+2\lambda_1\Vert\tilde{\phi}-u^*_N\Vert_{L^2(\Omega)}^2\\
    \leq&c_1(\Vert \tilde{\bpsi}-\nabla u^*_N\Vert_{H_{\divg}(\Omega)}^2+\Vert \tilde{\bpsi}-\nabla u^*_N\Vert_{H^1(\Omega)}^2)+c_2\Vert \tilde{\phi}-u^*_N\Vert_{H^1(\Omega)}^2
\end{aligned}\end{equation*}
where $c_1=\max\{1,2\lambda_1,\lambda_2C_{\Omega}\}$ and $c_2=\max\{2,2\lambda_1\}$. Then we can take minimum on the right hand side of the inequality and finish the proof.
\end{proof}
Now, we can prove the error estimates in Theorem \ref{theorem:neumman} using Lemma \ref{lemma:posterror_1} and Lemma \ref{lemma:losserror_1} directly.
\begin{proof}
Suppose
$$(\phi^*,\bpsi^*)=\underset{\phi\in V,\bpsi\in W}{\arg\min} L^N(\phi,\bpsi).$$
By Lemma \ref{lemma:posterror_1}, we have
\begin{equation*}
\Vert\hat{\phi}-u^*_N\Vert^2_{H^1(\Omega)}\lesssim L^N(\hat{\phi},\hat{\bpsi}),
\quad \text{and} \quad
\Vert\hat{\bpsi}-\nabla u^*_N\Vert_{H_{\divg}(\Omega)}^2\lesssim L^N(\hat{\phi},\hat{\bpsi}).
\end{equation*}
Since
\begin{equation*}\begin{aligned}
    L^N(\hat{\phi},\hat{\bpsi})&=L^N(\hat{\phi},\hat{\bpsi})-L_n^N(\hat{\phi},\hat{\bpsi})+L_n^N(\hat{\phi},\hat{\bpsi})\\
    &\leq L^N(\hat{\phi},\hat{\bpsi})-L_n^N(\hat{\phi},\hat{\bpsi})+L_n^N(\phi^*,\bpsi^*)\\
    &=L^N(\hat{\phi},\hat{\bpsi})-L_n^N(\hat{\phi},\hat{\bpsi})+L_n^N(\phi^*,\bpsi^*)-L^N(\phi^*,\bpsi^*)+L^N(\phi^*,\bpsi^*)\\
    &\leq 2\sup_{\phi\in V,\bpsi\in W}|L^N(\phi,\bpsi)-L^N_n(\phi,\bpsi)|+L^N(\phi^*,\bpsi^*)
    .
\end{aligned}\end{equation*}
Then we just apply Lemma \ref{lemma:losserror_1} and we finish the proof. 
\end{proof}
\subsection{Dirichlet boundary condition}
Problem (\ref{eq:pde_dirichlet}) also has a similar error estimate.
\begin{theorem}\label{theorem:dirchlet}
Suppose the solution $u^*_D$ to the Dirichlet problem (\ref{eq:pde_dirichlet}) and
% $$(\phi^*,\bpsi^*)=\underset{\phi\in V,\bpsi\in W}{\arg\min} L^D(\phi,\bpsi),$$
$$(\hat{\phi},\hat{\bpsi})=\underset{\phi\in V,\bpsi\in W}{\arg\min} L^D_n(\phi,\bpsi).$$
Assume parameters of networks in $W$ are uniformly bounded.
Then, we have
\begin{equation*}
    \Vert\hat{\phi}-u^*_D\Vert^2_{H^1(\Omega)}\lesssim\sqrt{E_3+E_4},
\quad \text{and} \quad
    \Vert\hat{\bpsi}-\nabla u^*_D\Vert_{H_{\divg}(\Omega)}^2\lesssim\sqrt{E_3+E_4},
\end{equation*}
where
$$E_3=2\sup_{\phi\in V,\bpsi\in W}|L^D(\phi,\bpsi)-L^D_n(\phi,\bpsi)|,$$
and
$$E_4=\underset{\bpsi\in W}{\min}\Vert \bpsi-\nabla u^*_D\Vert^2_{H_{\divg}(\Omega)}+\underset{\phi\in V}{\min}\Vert\phi-u^*_D\Vert^2_{H^1(\Omega)}.$$

\end{theorem}
Similar to the proof of Theorem \ref{theorem:neumman}, we also need to show the error of neural networks can be bounded by the expected loss.
\begin{lemma}\label{lemma:posterror_2}
Let $\phi\in V$ and $\bpsi\in W$ be two networks, and $u^*_D$ be the classical solution to the Dirichlet problem (\ref{eq:pde_dirichlet}). Then
\begin{equation*}
\Vert\phi-u^*_D\Vert^2_{H^1(\Omega)}\leq 2\tilde{C}_{\Omega}\left(C_BR_b^D+(R^D_b)^2+\frac{\tilde{C}_{\Omega}}{2}R_g^2+\frac{\tilde{C}_{\Omega}}{2}(R^D_{e})^2\right),
\end{equation*}
and
\begin{equation*}
\Vert\psi-\nabla u^*_D\Vert_{H_{\divg}(\Omega)}^2\leq 4C_B\tilde{C}_{\Omega}R^D_b+4\tilde{C}_{\Omega}(R^D_b)^2+(2+2\tilde{C}_{\Omega}^2)R^2_g+2\tilde{C}_{\Omega}^2(R^D_{e})^2+(R_{e}^D)^2,
\end{equation*}
where $C_B=\sqrt{|\partial\Omega|}\sup_{\bpsi\in W,x\in\partial\Omega}|\bpsi(x)|+C_{\Omega}\Vert\nabla u^*_D\Vert_{H^{1}(\Omega)}$ is a constant. 
\end{lemma}
\begin{proof}
We first need to introduce the Poincaré-Friedrichs inequality:
For any $u\in H^1(\Omega)$, there exist a constant $\tilde{C}_{\Omega}$ such that 
\begin{equation*}
    \Vert u\Vert^2_{H^1}\leq \tilde{C}_{\Omega}\left(\int_{\partial\Omega}|u|^2\ ds+\int_{\Omega}|\nabla u|^2\ dx\right).
\end{equation*}
Since $u^*_D$ is the classical solution to the Dirichlet problem (\ref{eq:pde_dirichlet}), we multiply $\hat{\phi}$ on both sides of (\ref{eq:pde_dirichlet}) and integrate them over $\Omega$, we get
\begin{equation*}
    -\int_{\partial\Omega}\hat{\phi}\nabla u^*_D\cdot\mathbf{n}\ ds+\int_{\Omega}\nabla u^*_D\nabla\hat{\phi}\ dx=\int_{\Omega}f\hat{\phi}\ dx.
\end{equation*}
By adding $\int_{\Omega}\divg(\bpsi)\hat{\phi}\ dx$ on both sides,
\begin{equation*}\begin{aligned}
    &\int_{\Omega}(f+\divg(\bpsi))\hat{\phi}\ dx\\
    =& -\int_{\partial\Omega}\hat{\phi}\nabla u^*_D\cdot\mathbf{n}\ ds+\int_{\Omega}\nabla u^*_D\nabla\hat{\phi}+\divg(\bpsi)\hat{\phi}\ dx\\
    =&-\int_{\partial\Omega}\hat{\phi}(\nabla u^*_D-\bpsi)\cdot\mathbf{n}\ ds+\int_{\Omega}(\nabla u^*_D-\bpsi)\nabla\hat{\phi}\ dx\\
    =&-\int_{\partial\Omega}(\phi-g_2)(\nabla u^*_D-\bpsi)\cdot\mathbf{n}\ ds+\int_{\Omega}(\nabla u^*_D-\nabla\phi+\nabla\phi-\bpsi)\nabla\hat{\phi}\ dx\\
    =&\int_{\partial\Omega}(\phi-g_2)(\bpsi-\nabla u^*_D)\cdot\mathbf{n}\ ds+\int_{\Omega}(\nabla u^*_D-\nabla\phi+\nabla\phi-\bpsi)\nabla\hat{\phi}\ dx.
\end{aligned}\end{equation*}
By rearranging the equation, we get
\begin{equation*}\begin{aligned}
    &\int_{\Omega}|\nabla\hat{\phi}|^2\ dx\\
    =&\int_{\partial\Omega}(\phi-g_2)(\bpsi-\nabla u^*_D)\cdot\mathbf{n}\ ds+\int_{\Omega}(\nabla\phi-\bpsi)\nabla\hat{\phi}\ dx-\int_{\Omega}(f+\divg(\bpsi))\hat{\phi}\ dx\\
    \leq&\Vert\bpsi-\nabla u^*_D\Vert_{L^2(\partial\Omega)}R_b^D+\frac{\tilde{C}_{\Omega}}{2}R_g^2+\frac{\Vert\nabla\hat{\phi}\Vert_{L^2(\Omega)}^2}{2\tilde{C}_{\Omega}}+\frac{\tilde{C}_{\Omega}}{2}(R^D_{e})^2+\frac{\Vert\hat{\phi}\Vert_{L^2(\Omega)}^2}{2\tilde{C}_{\Omega}}\\
    \leq&C_BR_b^D+\frac{\tilde{C}_{\Omega}}{2}R_g^2+\frac{\Vert\nabla\hat{\phi}\Vert_{L^2(\Omega)}^2}{2\tilde{C}_{\Omega}}+\frac{\tilde{C}_{\Omega}}{2}(R^D_{e})^2+\frac{\Vert\hat{\phi}\Vert_{L^2(\Omega)}^2}{2\tilde{C}_{\Omega}},
\end{aligned}\end{equation*}
where $C_B=\sqrt{|\partial\Omega|}\sup_{\bpsi\in W,x\in\partial\Omega}|\bpsi(x)|+C_{\Omega}\Vert\nabla u^*_D\Vert_{H^{1}(\Omega)}$. 

Using the Poincaré-Friedrichs inequality,
\begin{equation*}
    \begin{aligned}
        \Vert \hat{\phi}\Vert_{H^1(\Omega)}&\leq \tilde{C}_{\Omega}\left((R^D_b)^2+\Vert\nabla\hat{\phi}\Vert_{L^2(\Omega)}^2\right)\\
        &\leq \tilde{C}_{\Omega}\left((R^D_b)^2+C_BR_b^D+\frac{\tilde{C}_{\Omega}}{2}R_g^2+\frac{\tilde{C}_{\Omega}}{2}(R^D_{e})^2\right)+\frac{1}{2}\Vert\hat{\phi}\Vert_{H^1(\Omega)}^2,
    \end{aligned}
\end{equation*}
which implies that
\begin{equation*}
    \Vert\hat{\phi}\Vert_{H^1(\Omega)}^2\leq 2\tilde{C}_{\Omega}\left((R^D_b)^2+C_BR_b^D+\frac{\tilde{C}_{\Omega}}{2}R_g^2+\frac{\tilde{C}_{\Omega}}{2}(R^D_{e})^2\right),
\end{equation*}
and
\begin{equation*}
    \begin{aligned}
        \Vert\hat{\bpsi}\Vert_{H_{\divg}(\Omega)}^2&\leq \left(R_g+\Vert\nabla\hat{\phi}\Vert_{L^2(\Omega)}\right)^2+(R_{e}^D)^2\\
        &\leq 2R_g^2+(R_{e}^D)^2+2\Vert\hat{\phi}\Vert_{H^1(\Omega)}^2\\
        %&\leq 2R_g^2+(R_{e}^D)^2+4\tilde{C}_{\Omega}((R^D_b)^2+C_BR_b^D+\frac{\tilde{C}_{\Omega}}{2}R_g^2+\frac{\tilde{C}_{\Omega}}{2}(R^D_{e})^2)\\
        &\leq 4C_B\tilde{C}_{\Omega}R^D_b+4\tilde{C}_{\Omega}(R^D_b)^2+(2+2\tilde{C}_{\Omega}^2)R^2_g+2\tilde{C}_{\Omega}^2(R^D_{e})^2+(R^D_{e})^2.
    \end{aligned}
\end{equation*}
\ 
\end{proof}
Next, we can show the expected loss can be bounded by the approximation error.
\begin{lemma}\label{lemma:losserror_2}
Let
$$(\phi^*,\bpsi^*):=\underset{\phi\in V,\bpsi\in W}{\arg\min} L^D(\phi,\bpsi).$$
Then,
\begin{equation*}
    L^D(\phi^*,\bpsi^*)\lesssim \underset{\bpsi\in W}{\min}\Vert \bpsi-\nabla u^*_D\Vert^2_{H_{\divg}(\Omega)}+\underset{\phi\in V}{\min}\Vert\phi-u^*_D\Vert^2_{H^1(\Omega)}.
\end{equation*}
\end{lemma}
\begin{proof}
For any $\tilde{\phi}\in V$ and $\tilde{\bpsi}\in W$, we have
\begin{equation*}\begin{aligned}
    L^D(\phi^*,\bpsi^*)\leq& L^D(\tilde{\phi},\tilde{\bpsi})\\
    =&\Vert \tilde{\bpsi}-\nabla\tilde{\phi}\Vert_{L^2(\Omega)}^2+\lambda_1\Vert\divg(\tilde{\bpsi})+f\Vert_{L^2(\Omega)}^2+\lambda_2\Vert \tilde{\phi}-g_2\Vert_{L^2(\partial\Omega)}^2\\
    \leq&\lambda_1\Vert\divg(\tilde{\bpsi})-\divg(\nabla u^*_D)\Vert_{L^2(\Omega)}^2+\lambda_2C_{\Omega}\Vert\tilde{\phi}-u^*_D\Vert_{H^1(\Omega)}^2\\
    &+2\Vert \tilde{\bpsi}-\nabla u^*_D\Vert^2_{L^2(\Omega)}+2\Vert \nabla\tilde{\phi}-\nabla u^*_D\Vert^2_{L^2(\Omega)}\\
    \leq& C_3\Vert \tilde{\bpsi}-\nabla u^*_D\Vert^2_{H_{\divg}(\Omega)}+C_4\Vert \tilde{\phi}-u^*_D\Vert^2_{H^1(\Omega)},
\end{aligned}\end{equation*}
where $C_3=\max\{2,\lambda_1\}$ and $C_4=\max\{2,\lambda_2C_{\Omega}\}$. 
\end{proof}
Then, the proof of Theorem \ref{theorem:dirchlet} is quite similar to Theorem \ref{theorem:neumman}. It follows directly from Lemma \ref{lemma:posterror_2} and \ref{lemma:losserror_2}.

\subsection{Convergence Rate and Comparison with DRM}
Suppose the the approximation error of neural networks in the set $V$ to any function $u$ in a specific space $X$, equipped with the norm $\Vert\cdot\Vert_X$, has the following estimates
\begin{equation*}
    \min_{\phi\in V}\Vert \phi-u\Vert_{H^1(\Omega)}^2\lesssim m^{-\gamma}\Vert u\Vert_X^2,
\end{equation*}
where $m$ is a parameter proportional to the total number of parameters in networks of $V$ and $\gamma>0$ is a positive constant. Similarly, we can expect the approximation error of networks in $W$ to $\nabla u$ has the same rate if $\nabla u$ has desired regularity:
\begin{equation*}
    \min_{\bpsi\in W}\Vert \bpsi-\nabla u\Vert_{H^1(\Omega)}^2\lesssim m^{-\gamma}\Vert \nabla u\Vert_{X^d}^2,
\end{equation*}
where $\Vert \nabla u\Vert_{X^d}^2=\sum_{i=1}^d \Vert\partial_i u\Vert_X^2.$ Besides, $\Vert \bpsi-\nabla u\Vert_{H_{\divg}(\Omega)}^2\leq d\Vert\bpsi-\nabla u\Vert_{H^1(\Omega)}^2$, so we have
\begin{equation*}
    \min_{\bpsi\in W}\Vert \bpsi-\nabla u\Vert_{H_{\divg}(\Omega)}^2\lesssim \frac{1}{m^\gamma}\Vert \nabla u\Vert_{X^d}^2.
\end{equation*}
Then we can derive the convergence rate for the mixed residual method when solving the Neumann problem (\ref{eq:elliptic_N}) and the Dirichlet problem (\ref{eq:elliptic_D}). Let us neglect the quadrature error part here, since the estimates of quadrature errors are quite similar for different methods. For the Neumann problem, let $(\phi^*,\bpsi^*):=\underset{\phi\in V,\bpsi\in W}{\arg\min} L^N(\phi,\bpsi)$ and $u^*_N\in X$. The following convergence rate can be derived from Theorem \ref{theorem:neumman}
$$\Vert\phi^*-u^*_N\Vert_{H^1(\Omega)}^2\lesssim\frac{1}{m^\gamma} \text{ and } \Vert\bpsi^*-\nabla u^*_N\Vert_{H_{\divg}(\Omega)}^2\lesssim\frac{1}{m^\gamma}.$$
For the Dirichlet problem, by abusing the notation, let $(\phi^*,\bpsi^*):=\underset{\phi\in V,\bpsi\in W}{\arg\min} L^D(\phi,\bpsi)$ and $u^*_D\in X$. Then, from Theorem \ref{theorem:dirchlet}, we can derive the convergence rate:
$$\Vert\phi^*-u^*_D\Vert_{H^1(\Omega)}^2\lesssim\frac{1}{m^{\gamma/2}} \text{ and } \Vert\bpsi^*-\nabla u^*_D\Vert_{H_{\divg}(\Omega)}^2\lesssim\frac{1}{m^{\gamma/2}}.$$

Let us compare the above convergence rate with DRM. The priori error of DRM with homogeneous boundary conditions has been well-studied in the literature \cite{lu2021priori,li2021generalization,muller2021error}. 
For the Neumann problem with homogeneous boundary condition, 
%DRM uses the following loss function:
% \begin{equation*}
%     J(u)=\Vert u-f\Vert_{L^2(\Omega)}^2+\Vert\nabla u\Vert_{L^2(\Omega)}^2.
% \end{equation*}
let $\phi^*$ be DRM solution. Then, the error of $\phi^*$ can be estimated by
\begin{equation*}
    \Vert\phi^*-u^*_N\Vert_{H^1(\Omega)}^2\leq\min_{\phi}\Vert\phi-u^*_N\Vert_{H^1(\Omega)}^2\lesssim\frac{1}{m^\gamma}\Vert u\Vert_X^2.
\end{equation*}
% If the network $\phi$ is a two-layer network in $V^m(u^*)$, then we have 
% $$\Vert\phi^*-u^*\Vert^2_{H^1}\lesssim \frac{\Vert u^*\Vert_{\B}^2}{m}.$$ 
For the Dirichlet problem with homogeneous boundary condition, 
%DRM works with the following loss:
% \begin{equation*}
%     J_\lambda(\phi)=-\langle\phi,f\rangle_{L^2(\Omega)}+\Vert\nabla \phi\Vert_{L^2(\Omega)}^2+\lambda\Vert u\Vert_{L^2(\partial\Omega)},
% \end{equation*}
% where $\lambda$ controls the weight of boundary penalty. 
let $\phi^*_\lambda$ be the DRM solution where $\lambda$ is the weight of boundary penalty. Then it has the following estimates in \cite{muller2021error}:
\begin{equation*}\begin{aligned}\Vert\phi^*_\lambda-u^*_D\Vert_{H^1(\Omega)}^2\lesssim\left(\sqrt{\frac{1+\lambda}{m^\gamma}\Vert u^*_D\Vert_X^2+\frac{1}{\lambda}}+\frac{1}{\lambda}\right)^2.\end{aligned}\end{equation*}
The best rate is then $\frac{1}{m^{\gamma/2}}$ when choosing $\lambda=m^{\gamma/2}$.

Generally, the mixed residual method achieve the same convergence rate with respect to $m$ with DRM in both problems. Our result suggests that the mixed residual method can also recover the Laplacian of the solution by taking divergence of the $\bpsi$ network. Moreover, our result can be applied to both homogeneous boundary and non-homogeneous boundary problems.

\section{Estimation of the approximation error and the quadrature error} \label{sec:barron}
In this section, we are going to estimate the approximation error and quadrature error respectively.
\subsection{Approximation error of Barron space functions}
By abusing the notation a little, let $V^m_u$ and $W^m_{\boldsymbol{u}}$ be the sets of two layer neural networks of width $m$ with respect to function $u$ or $\boldsymbol{u}$. The approximation error of networks in $V^m_u$ can be generalized by the result of \cite{li2021generalization}:
\begin{lemma}\label{theorem:H1sp}
For any function $u\in\B(\Omega)$ and $m\in\mathbb{N}^+$, there exist a network $u_m\in V^m_u$ such that
\begin{equation*}
    \Vert u-u_m\Vert_{H^1(\Omega)}^2\lesssim  \frac{\|u\|_{\B(\Omega)}^2}{m}.
\end{equation*}
\end{lemma}
It is easy to approximate a ReLU activated network by a ReQU activated network while it is not trivial to show that the new coefficients in the neural network is bounded. For this reason, it requires $u\in \mathcal{B}(\Omega)$. The detailed proof of this lemma is shown in appendix \ref{app_approx_er}.

Further, we can easily generalize this approximation error to functions in $\B(\Omega;\mathbb{R}^d)$ with the set of networks $W^m_{\boldsymbol{u}}$.
For any $\bpsi\in W^m_{\boldsymbol{u}}$, the $i$th entry of $\bpsi$ is a two-layer network in $V^m_{\boldsymbol{u}_i}$ and the following lemma is a direct corollary of lemma \ref{theorem:H1sp}.
\begin{lemma}\label{theorem:H1spV}
For any function $\boldsymbol{u}\in\B(\Omega;\mathbb{R}^d)$ and $m\in\mathbb{N}^+$, there exist a network $\boldsymbol{u_m}\in W^m_{\boldsymbol{u}}$ such that
\begin{equation*}
    \Vert \boldsymbol{u}-\boldsymbol{u_m}\Vert^2_{H^1(\Omega)}\lesssim \frac{\Vert \boldsymbol{u}\Vert^2_{\B(\Omega;\mathbb{R}^d)}}{m}.
\end{equation*}
\end{lemma}

\subsection{Estimation of the Rademacher Complexity}{\label{sec:rc_complexity}}
In this section, we will illustrate that the quadrature error trained on a finite dataset $\{X_i\}_{i=1}^n$ can be estimated by the Rademacher complexity. For simplicity, we consider $X_i\in\Omega = [0,1]^d$ in the following estimation.
\begin{definition}
    Let $\{X_i\}_{i=1}^n$ be a set of random variables independently distributed and $\{\varepsilon_i\}_{i=1}^n$ be an i.i.d sequence of Rademacher variables (i.e. taking the values of \{1,-1\} equiprobably). Then the \textbf{empirical Rademacher Complexity} of the function class $\mathcal{F}$ is a random variable given by 
    \begin{equation*}
      \hat{R}_n(\mathcal{F}):= \mathbb{E}_\varepsilon[\sup_{f\in\mathcal{F}}|\frac{1}{n}\sum_{i=1}^n \varepsilon_i f(X_i)|].
    \end{equation*}
    Taking its expectation yields the \textbf{Rademacher Complexity} of the function class $\mathcal{F}$
    \begin{equation*}
      R_n(\mathcal{F}):= \mathbb{E}_X[\hat{R}_n(\mathcal{F})]=\mathbb{E}_X\mathbb{E}_\varepsilon[\sup_{f\in\mathcal{F}}|\frac{1}{n}\sum_{i=1}^n \varepsilon_i f(X_i)|].
    \end{equation*}
\end{definition}

The definition gives some basic calculation rules:
\begin{lemma}
    Let $\mathcal{F},\mathcal{G}$ be function classes and $a,b$ be constants. Then
    \begin{itemize}
        \item [(i)] $R_n(\mathcal{F+G})\leq R_n(\mathcal{F}) + R_n(\mathcal{G})$.
        \item [(ii)] $R_n(a\mathcal{F}) = |a|R_n(\mathcal{F}).$
        \item [(iii)] Assume $g$ is a fixed function and $\|g\|_\infty\leq b$, then $R_n(g)\leq \frac{b}{\sqrt{n}}.$
        \item [(iv)](Ledoux-Talagrand contraction lemma\cite{ledoux1991probability}) Assume that $\sigma:\mathbb{R}\mapsto\mathbb{R}$ is $l$-Lipschitz with $\sigma(0)=0$, then $R_n(\sigma(\mathcal{F}))\leq 2l R_n(\mathcal{F})$.
        \item [(v)] $R_n(\mathcal{F}^2)\leq 4\sup_{f\in\mathcal{F}}\|f\|_\infty R_n(\mathcal{F})$.
        \item [(vi)] $R_n(\mathcal{F}\mathcal{G})\leq 6\sup_{f\in\mathcal{F}\cup\mathcal{G}}\|f\|_{\infty}(R_n(\mathcal{F})+R_n(\mathcal{G}))$.
    \end{itemize}
    \label{lem::computing_rules}
\end{lemma}
The proof can be found in Appendix \ref{app_CR}.
With the calculation rules prepared, we are ready to estimate the complexity of the neural network function classes and the loss function classes.

\begin{lemma}
    Let $\mathcal{G}$ be the linear transformation function class defined by 
    \begin{equation*}
        \mathcal{G}:= \{\omega\cdot x + b|\|\omega\|_2=1, |b| \leq 1\}. 
    \end{equation*}
    Then we have 
    \begin{equation}
        R_n(\mathcal{G})\leq \frac{\sqrt{2d\log d}+1}{\sqrt{n}}.
    \end{equation}
    \label{lem::linear_trans}
\end{lemma}
\begin{proof}
    It is obvious that the Rademacher complexity of the constant is bounded by $\frac{1}{\sqrt{n}}$. Then we consider the function class 
    \begin{equation*}
        \hat{\mathcal{G}}:=\{\omega\cdot x|\|w\|_2=1\}.
    \end{equation*}
    Dividing the Rademacher complexity into each dimensional yields
    \begin{equation*}
        %\begin{aligned}
            R_n(\hat{\mathcal{G}})\leq \frac{\|w\|_1}{n}\mathbb{E}_X\mathbb{E}_\varepsilon\|\sum_{i=1}^n\varepsilon_iX_i\|_{\infty}\leq \frac{\sqrt{d}\|w\|_2}{n}\mathbb{E}_X\mathbb{E}_\varepsilon\|\sum_{i=1}^n\varepsilon_iX_i\|_{\infty}.
        %\end{aligned}
    \end{equation*}
    Then we use the Massart lemma\cite{massart2000some}:
    \begin{equation}
        \mathbb{E}_\varepsilon[\max_{a\in A}\frac{1}{n}\sum_{i=1}^n\varepsilon_i a_i]\leq \max_{a\in A} \|a\|_2 \frac{\sqrt{2\log(|A|)}}{n}.
        \label{lem::massart}
    \end{equation}
    View each fearture $X_{i,j}$ and $-X_{i,j}$ for $j=1,\cdots,d$ as a member of the finite hypothesis class, i.e., 
    \begin{equation*}
        \begin{aligned}
            A &= \{Y_1,\cdots, Y_d,-Y_1,\cdots, -Y_d\},\\
            Y_j &= (X_{1,j},\cdots, X_{n,j}), \qquad j = 1,\cdots,d.
        \end{aligned}
    \end{equation*} 
    Applying (\ref{lem::massart}) yields
    \begin{equation*}
        \begin{aligned}
            &\mathbb{E}_X\mathbb{E}_\varepsilon\|\frac{1}{n}\sum_{i=1}^n\varepsilon_iX_i\|_\infty=\frac{1}{n}\mathbb{E}_X\mathbb{E}_\varepsilon\max_{a\in A}\sum_{i=1}^n\varepsilon_ia_i\\
            &\leq \frac{\sqrt{2\log 2d}}{n}\mathbb{E}_X\max_{a\in A} \|a\|_2 \leq \frac{\sqrt{2\log 2d}}{\sqrt{n}}.
        \end{aligned}
        \label{eq::EX}
    \end{equation*}
    The last step is finished by $\|Y_i\|_2\leq \sqrt{n}\|Y_i\|_\infty = \sqrt{n}$, for all $j = 1,\cdots,d$. Summarizing the result, we have the inequality
    \begin{equation*}
        R_n(\mathcal{G})\leq \frac{\sqrt{2d\log 2d}+1}{\sqrt{n}}.
    \end{equation*}
    \ 
\end{proof}

The estimation of the complexity of a two-layer neural network depends on the activation function, so we make the following assumptions: $\sigma\in H^2(\Omega)$ with $\sigma(0)=0$, and
\begin{equation*}
         \sup_x |\sigma^{(k)}(x)|\leq l_k, \ k=0,1,2.
\end{equation*}

% With these properties, we are able to estimate the Rademacher complexity of complex loss functions. Before we start, the following assumptions are given:
% \begin{itemize}
%     \item[(i)] The activation function $\sigma\in C^2(\Omega)$ satisfies $\sigma(0)=0$ and \begin{equation*}
%         \sup_x |\sigma^{(k)}(x)|\leq l_k, \ k=0,1,2.
%     \end{equation*}
%     \item[(ii)] $\sup_x|f(x)|\leq \normu.$
%     \item[(iii)] $\sup_x|g(x)|\leq \normu.$ 
% \end{itemize}

\begin{lemma}
    The Rademacher complexity of $V^m_{u}$ is bounded by 
    \begin{equation}
        R_n(V^m_u) \leq \frac{C_1\normu}{\sqrt{n}},
    \end{equation} 
    where $C_1$ depends on $d,$ $l_1$.
    \label{lem::RC_NN}
\end{lemma}
\begin{proof}
    Using properties (i) and (iv), the Rademacher complexity of two layer neural network is broken down into the sum of the Rademacher complexity of each neuron, that is  
    \begin{equation*}
        \begin{aligned}
            R_n(V^m_u) &\leq \frac{2\normu}{\sqrt{n}} + \frac{1}{m}\sum_{i=1}^m |a_i|R_n(\sigma(\mathcal{G}))
        \end{aligned}
    \end{equation*}
    With the assumption that $\sigma(x)$ is $l_1$-Lipschitz and $\sigma(0)=0$, we have 
    \begin{equation*}
        \begin{aligned}
            R_n(V^m_u)& \leq \frac{2\normu}{\sqrt{n}} + \frac{1}{m}\sum_{i=1}^m |a_i|2l_1R_n(\mathcal{G}_i)\\
            &\leq \frac{2\normu}{\sqrt{n}} + \frac{2l_1}{m}\sum_{i=1}^m |a_i|\sup_i R_n(\mathcal{G}_i)\\
            &\leq \frac{2\normu}{\sqrt{n}} + 16 l_1 R_n(\mathcal{G})\normu.\\
        \end{aligned}
    \end{equation*}
    Using the result of lemma \ref{lem::linear_trans} concludes that 
    \begin{equation*}
        R_n(V^m_u)\leq \frac{(16l_1+2 + 16l_1\sqrt{2d\log 2d})\normu}{\sqrt{n}}.
    \end{equation*}
    \ 
\end{proof}

Then we further assume that $\|f\|_{L^\infty(\Omega)}\leq \|u\|_{\B^4(\Omega)}$ and $\|g\|_{L^\infty(\partial\Omega)}\leq \|u\|_{\B^4(\Omega)}$.
In order to complete the quadrature error, we give the estimation of Rademacher complexity with respect to the following function classes which are determined by the expected loss functions $L^N(\phi, \psi)$ and $L^D(\phi,\psi)$. 
\begin{equation}
    \begin{aligned}
        \mathcal{L}^N&:=\{|\nabla \phi-\psi|^2 + \lambda_1(-\nabla\cdot\psi+\phi-f)^2\vert\phi\in V^m_u,\psi\in W^m_{\nabla u}\},\\
        \mathcal{L}_b^N& = \{\lambda_2(\psi\cdot \mathbf{n}-g_1)^2\vert\psi\in W^m_{\nabla u}\},\\
    \mathcal{L}^D&:=\{|\nabla \phi-\psi|^2 + \lambda_1(-\nabla\cdot\psi-f)^2 \vert\phi\in V^m_u,\psi\in W^m_{\nabla u}\},\\
    \mathcal{L}_b^D& = \{\lambda_2(\phi-g_2)^2\vert\phi\in V^m_u\}.\\
\end{aligned}
\label{eq::loss_function_class}
\end{equation}
When $u\in\B^4(\Omega)$, $\nabla u \in \B(\Omega;\mathbb{R}^d)$ and lemma \ref{lem::vector_bnorm} implies $\|\nabla u\|_{\B(\Omega;\mathbb{R}^d)}\leq \|u\|_{\B^4(\Omega)}$. 
\begin{lemma}[Rademacher complexity for the Neumann and the Dirichlet problem]\label{lem::rc_of_neumann}
    Let $\mathcal{L}^N,\mathcal{L}^N_b,\mathcal{L}^D,\mathcal{L}^D_b$ defined as (\ref{eq::loss_function_class}) be the function classes of deep mixed residual method with respect to the two layer networks $V_u^m$ and $W_{\nabla u}^m$. Then 
    \begin{equation}
        \begin{aligned}
            R_n(\mathcal{L}^N) + R_{\bar{n}}(\mathcal{L}^N_b) \leq \frac{C_2\|u\|_{\B^4(\Omega)}^2}{\sqrt{n}},\\
            R_n(\mathcal{L}^D) + R_{\bar{n}}(\mathcal{L}^D_b) \leq \frac{C_2\|u\|_{\B^4(\Omega)}^2}{\sqrt{n}},
        \end{aligned}
    \end{equation}
    where $C_2$ depends on $d$, $l_0$, $l_1$, $l_2$ and $\bar{n} \sim \frac{n}{d^2}$.
\end{lemma}
\begin{proof}
    According to the definition of the Rademacher complexity, it follows directly that the upper bound of $L^N$ and $L^D$ can be dividide into the following parts:
    \begin{itemize}
        \item [1.] $\mathcal{L}_g=\{|\nabla \phi - \psi|^2|\phi\in V^m_u,\psi\in W^m_{\nabla u}\}$,       
        \item [2.] $\mathcal{L}^N_{e}=\{(-\nabla\cdot\psi+\phi-f)^2|\phi\in V^m_u,\psi\in W^m_{\nabla u}\}$,
        \item [3.] $\mathcal{L}^D_{e}=\{((-\nabla\cdot\psi-f)^2|\psi\in W^m_{\nabla u}\} $.
    \end{itemize}
    The first term appears as a gradient penalty term in the loss function of both problems. By expanding the first term, we obtain
%    \begin{equation*}
        $|\nabla\phi - \psi|^2 = \sum_{j=1}^d (\phi_{x_j}+\psi_j)^2.$
%    \end{equation*}
    Viewing each dimension yields 
    \begin{equation*}
        %\begin{aligned}
            \phi_{x_j} = \frac{1}{m}\sum_{i=1}^na_iw_{i,j}\sigma'(w_i\cdot x + b_i), \ \psi_j\in V^m_u, \ j = 1,\cdots,d.
        %\end{aligned}
    \end{equation*}
    It follows that 
    \begin{equation*}
            R_n(\mathcal{L}_g)\leq 4\sum_{j=1}^d \sup_{\phi\in V^m_u,\psi\in W^m_{\nabla u} } \Vert \nabla \phi - \psi\Vert_{\infty}(R_n(V^m_u)+\frac{1}{m}\sum_{i=1}^m|a_iw_{i,j}|R_n(\sigma'(\mathcal{G}))), 
    \end{equation*}
    where 
    \begin{equation*}
        \begin{aligned}
            \sup_{\phi\in V^m_u,\psi\in W^m_{\nabla u}} \Vert \nabla \phi - \psi\Vert_{\infty}&\leq 
            (8l_1+ 2+ 8l_0)\|u\|_{\B^4(\Omega)},\\
            \frac{1}{m}\sum_{j=1}^d\sum_{i=1}^m |a_iw_{i,j}|R_n(\sigma'(\mathcal{G}))&\leq \frac{1}{m}\sum_{i=1}^m |a_i| \sqrt{d}\|w_i\|_2R_n(\sigma'(\mathcal{G}))\\
            &\leq 8\sqrt{d}\|u\|_{\B^4(\Omega)} R_n(\mathcal{\sigma'(G)}) \\
            &\leq 16\sqrt{d}l_2\|u\|_{\B^4(\Omega)} R_n(\mathcal{G}).
        \end{aligned}
    \end{equation*}
    According to lemma \ref{lem::linear_trans} and \ref{lem::RC_NN}, it follows that 
    \begin{equation}
        R_n(\mathcal{L}_g)\lesssim \sqrt{\frac{d^3\log d}{n}}\|u\|_{\B^4(\Omega)}^2.
        \label{eq::gradient_penalty}
    \end{equation} 
    The estimation of the second part follows in a similar manner.
   \begin{equation*}
       \begin{aligned}
           R_n(\mathcal{L}_{e}^N)&\leq 4 \sup_{\phi\in V^m_u,\psi\in W^m_{\nabla u} }\Vert -\nabla\cdot \psi+\phi-f \Vert_\infty\\
           &\cdot(8d\|u\|_{\B^4(\Omega)} R_n(\mathcal{\sigma'(G)})+R_n(V_u^m)+R_n(f))
       \end{aligned}
   \end{equation*}  
   By the fact that $\psi_i \in V_u^m$, the estimation of $\nabla\cdot\psi$ is almost equivalent to that of $\phi_{x_i}$, which gives 
   \begin{equation*}
    \sup_{\phi\in V^m_u,\psi\in W^m_{\nabla u} }\Vert -\nabla\cdot \psi+\phi-f \Vert_\infty\leq (8dl_1 + 3 + 8l_0)\|u\|_{\B^4(\Omega)}.
   \end{equation*}
    Then the second part can be proved in the same way as shown before, and it is easy to show 
    \begin{equation}
        R_n(\mathcal{L}_{e}^N)\lesssim \sqrt{\frac{d^5\log d}{n}}\|u\|_{\B^4(\Omega)}^2.
        \label{eq::neumann_pde}
    \end{equation}
    $\mathcal{L}^D_{e}$ and $\mathcal{L}^N_{e}$ have almost the same form. Consequently, the complexity differs by a constant factor. 
    
    The last step is to estimate the function classes with respect to the boundary conditions. With the assumption that the domain $\Omega$ is a rectangle, it follows 
    \begin{equation}
        \begin{aligned}
            R_{\bar{n}}(\mathcal{L}_{b}^N)&\leq 4\sup_{\psi\in W^m_{\nabla u}} \Vert\psi\cdot \mathbf{n}-g_1\Vert_{L^{\infty}(\partial\Omega)}(R_{\bar{n}}(V^m_u)+R_{\bar{n}}(g_1))\lesssim  \sqrt{\frac{d\log d }{{\bar{n}}}}\|u\|_{\B^4(\Omega)}^2.
        \end{aligned}
        \label{eq::neumann_bc}
    \end{equation}
    In general, the number of sample points on the boundary is different from that in the interior domain. Comparing with (\ref{eq::gradient_penalty}) and (\ref{eq::neumann_pde}), taking $\bar{n}=\frac{n}{d^2}$ will not change the upper bound of the Rademacher complexity.  
    
     The formula (\ref{eq::gradient_penalty}), (\ref{eq::neumann_pde}) and (\ref{eq::neumann_bc}) conclude a upper bound of the Rademacher complexity.
    \begin{equation}
        R_n(\mathcal{L}^N) +R_{\bar{n}}(\mathcal{L}^N_b)\lesssim \sqrt{\frac{d^5\log d}{n}}\|u\|_{\B^4(\Omega)}^2.
    \end{equation}
    
    The estimation for the Dirichlet problem is almost identical. Although these are two different problems, they have the same upper bound, which is determined by Laplacian, dimensionality, and size of the dataset.
\end{proof}

    It remains to show that the Rademacher complexity can bound the quadrature error. 
    % For better understanding, a one-variable case is given. Let $L(u)$ be an arbitrary given loss function and $L_n(u) = \frac{1}{n}\sum_{i=1}^n l(u(X_i))$ be the corresponding empirical loss function under a dataset $\mathcal{D}=\{X_1,\cdots,X_n\}$. It satisfies $\mathbb{E} L_n(u)= L(u)$. For example, if $L(u) = \int_\Omega u^2(x) \mathrm{d}x$, then $L_n(u) = \frac{1}{n}\sum_{i=1}^n u(x_i)^2$ and $l(u) = u^2$.
    The following lemma fills up the gap between the Rademacher complexity and the quadrature error.
 \begin{lemma}\label{lemma:quadrature_error}
    Let $\mathcal{F}$ be a set of functions, $\{X_1,\cdots,X_n\}$ be i.i.d. random variables, and $\mathbb{E}_X \frac{1}{n}\sum_{i=1}^nl(u(X_i))= L(u)$. Then
     \begin{equation}
         \mathbb{E}_X\sup_{u\in\mathcal{F}}|L(u) - \frac{1}{n}\sum_{i=1}^nl(u(X_i))|\leq 2R_n(\mathcal{L}),
     \end{equation}
     where $\mathcal{L}:=\{l(u)|u\in\mathcal{F}\}$.
 \end{lemma}
 \begin{proof}
    We first show the inequality that 
    \begin{equation*}
      \sup_f \mathbb{E}_Xf(X)\leq  \mathbb{E}_X \sup_f f(X).  
    \end{equation*}
    It follows directly from 
    $\mathbb{E}_Xf(x)\leq \mathbb{E}_X\sup_f f(X) $ for all $f$. Recall that 
        $L(u) = \frac{1}{n}\mathbb{E}_X\sum_{i=1}^n l(u(X_i)),$
    then 
     \begin{equation*}
        \begin{aligned}
                    \mathbb{E}_X\sup_{u\in\mathcal{F}}|L(u) - \frac{1}{n}\sum_{i=1}^nl(u(X_i))|&=\mathbb{E}_X\sup_{u\in\mathcal{F}}|\mathbb{E}_Y\frac{1}{n}\sum_{i=1}^n (l(u(Y_i))-l(u(X_i)))|\\
                    &\leq \mathbb{E}_X\mathbb{E}_Y\sup_{u\in\mathcal{F}}|\frac{1}{n}\sum_{i=1}^n (l(u(Y_i))-l(u(X_i)))|.
        \end{aligned}
     \end{equation*}
     Here $X_i$ and $Y_i$ have the same distribution, which yields
     \begin{equation*}
        \begin{aligned}
          &\quad \mathbb{E}_X\mathbb{E}_Y\sup_{u\in\mathcal{F}}|\frac{1}{n}\sum_{i=1}^n (l(u(Y_i))-l(u(X_i)))|\\
        &=   \mathbb{E}_X\mathbb{E}_Y\mathbb{E}_\epsilon\sup_{u\in\mathcal{F}}|\frac{1}{n}\sum_{i=1}^n \epsilon_i(l(u(Y_i))-l(u(X_i)))|  \\
        &\leq \mathbb{E}_Y\mathbb{E}_\epsilon\sup_{u\in\mathcal{F}}|\frac{1}{n}\sum_{i=1}^n \epsilon_il(u(Y_i))| + \mathbb{E}_X\mathbb{E}_\epsilon\sup_{u\in\mathcal{F}}|\frac{1}{n}\sum_{i=1}^n \epsilon_il(u(X_i))|\\
        &= 2R_n(\mathcal{L}).
        \end{aligned}
     \end{equation*} 
     \  
 \end{proof}

\subsection{Proof of the main theorems}
Now we can use the previous results to prove the main theorems. The proof of Theorem \ref{theorem:main_N} is given below.
\begin{proof}
From Theorem \ref{theorem:neumman}, we know that the error $\Vert\hat{\phi}-u^*_N\Vert_{H^1(\Omega)}^2$ and $\Vert\hat{\bpsi}-\nabla u^*_N\Vert_{H_{div}(\Omega)}^2$ can be bounded by quadrature error:
$$E_1=2\sup_{\phi\in V,\bpsi\in W}|L^N(\phi,\bpsi)-L^N_n(\phi,\bpsi)|$$
and approximation error:
\begin{equation*}\begin{aligned}
&E_2=\underset{\bpsi\in W}{\min}\left(\Vert\bpsi-\nabla u^*_N\Vert_{H_{\divg}(\Omega)}^2+\Vert \bpsi-\nabla u^*_N\Vert_{H^1(\Omega)}^2\right)+\underset{\phi\in V}{\min}\left(\Vert \phi-u^*_N\Vert_{H^1(\Omega)}^2\right).
\end{aligned}\end{equation*}
By Lemma \ref{theorem:H1sp} and \ref{theorem:H1spV}, since $u^*_N\in\B^4(\Omega)$ and $\nabla u^*_N\in\B(\Omega;\mathbb{R}^d)$, we can derive an estimate for the approximation error: $E_2\lesssim \frac{\Vert u^*_N\Vert^2_{\B^4(\Omega)}}{m}.$
For the quadrature error, it is easy to check the activation function $\text{ReQU}(x)$ satisfies
\begin{equation*}
    \begin{aligned} 
        |\text{ReQU}(\omega\cdot x + b)| &\leq (\sqrt{d}+1)^2,\\
        |\text{ReQU}'(\omega\cdot x + b)| &\leq 2(\sqrt{d}+1),\\
        |\text{ReQU}''(\omega\cdot x +b)| &\leq 2,
    \end{aligned}
\end{equation*}
for all $x\in\Omega$. Then we can apply Lemma \ref{lem::rc_of_neumann} and \ref{lemma:quadrature_error} to derive the estimate hold for any $\phi\in V^m_{u^*_N}$ and $\bpsi\in W^m_{u^*_N}$:
$$\mathbb{E}|L^N_n(\phi,\bpsi)-L^N(\phi,\bpsi)|\leq 2(R_n(\mathcal{L}^N)+R_{\bar{n}}(\mathcal{L}^N_b))\lesssim \frac{\Vert u^*_N\Vert_{\B^4(\Omega)}^2}{\sqrt{n}},$$
if we take $\bar{n}$ no less than $n/d^2$ and the constant is at most a polynomial of $d$. Then, $E_1$ also has the same estimate. Therefore, we have proved the theorem.
\end{proof}
The Theorem \ref{theorem:main_D} can also be proved in the same way.

\section{Numerical experiment}\label{sec:experiment}
\subsection{Experiment setup}
In this part, we are going to compare DRM with the mixed residual methods for solving various PDEs. The DGM is used as a benchmark. We use ResNet for all the experiments. The activation function is the RePU function:
\begin{equation*}
    \text{RePU}(x):=(\max\{0,x\})^p,
\end{equation*}
where $p\geq 2$. It is also called ReQU and ReCU when $p=2$ and $p=3$ respectively. To control the total number of parameters, the depth of the networks is fixed to be 10 in all the experiments and the width depends on a given positive integer $w$. For the deep Ritz method and the deep Galerkin method, the width is set to be $\lceil \sqrt{5}w\rceil$. For the mixed residual method, the width of the $\phi$ network is $w$ and the width of the $\bpsi$ network is $2w$. For example, if we give $w=10$, then the network width is $23$ for DRM and DGM, and the width is $10$ and $20$ for the mixed residual method. 

In the training process, we randomly generate 1,000 samples from the interior of the domain and 1,000 samples from the boundary in each iteration. We then update the parameters using the Adam algorithm with a learning rate $10^{-4}$. For each test, the total number of iterations is 500,000. 

We evaluate the algorithms by computing the relative errors:
\begin{equation*}\begin{aligned}
    e_0=\sqrt{\frac{\int|\phi-u^*|^2 dx}{\int|u^*|^2dx}},\ e_1=\sqrt{\frac{\int|\nabla\phi-\nabla u^*|^2 dx}{\int|\nabla u^*|^2dx}},\ e_2=\sqrt{\frac{\int|\Delta\phi-\Delta u^*|^2 dx}{\int|\Delta u^*|^2dx}},
\end{aligned}\end{equation*}
where $u^*$ is the analytical solution to the original PDE problem and $\phi$ is the neural network solution. For the mixed residual method, $\nabla\phi$ and $\Delta\phi$ in $e_1$ and $e_2$ are replaced by $\bpsi$ and $\divg(\bpsi)$. The integral is estimated using 10,000 randomly sampled quadrature points.
\subsection{Elliptic PDE}
We consider the Poisson equation with Dirichlet boundary condition:
\begin{align}
    -\Delta u(x) = d\pi^2\prod_{i=1}^d\sin(\pi x_i) \label{eq:elliptic_D}\text{ in }\Omega, \quad u(x)= 0 \text{ on } \partial\Omega, \notag
\end{align}
where $\Omega=[0,1]^d$. The true solution is $u^*(x)=\prod_{i=1}^d\sin(\pi x_i)$. 
% The expected loss function for the mixed residual method is (\ref{eq:expected_loss_D}) where $\lambda_1=100$ and $\lambda_2=10$. The loss of DRM is
% \begin{equation*}
%     \Vert\nabla\phi\Vert_{L^2(\Omega)}^2-\langle \phi,d\pi^2\prod_{i=1}^d\sin(\pi x_i)\rangle_{L^2(\Omega)}+\lambda_1\Vert\phi-0\Vert^2_{L^2(\partial\Omega)},
% \end{equation*}
% where $\lambda_1=1000$, and the loss of DGM is
% \begin{equation*}
%     \Vert\Delta\phi+d\pi^2\prod_{i=1}^d\sin(\pi x_i)\Vert_{L^2(\Omega)}^2+\lambda_1\Vert\phi-0\Vert^2_{L^2(\partial\Omega)},
% \end{equation*}
% where $\lambda_1=100$. In practice, changing the weights of residual terms would not affect the results very much. 
The numerical result is shown in Table \ref{table:poisson_D}. We can observe that DRM has the largest relative error in all cases. Besides, the $e_2$ error of DRM is significantly larger than its $e_0$ and $e_1$ error. The mixed residual method achieves similar accuracy to the benchmark while it cost about $30\%$ less time. 

\begin{table}
\centering
\begin{tabular}{|p{0.08\textwidth}|p{0.07\textwidth}|p{0.1\textwidth}|p{0.1\textwidth}|p{0.1\textwidth}|p{0.1\textwidth}|p{0.09\textwidth}|p{0.08\textwidth}|} 
\hline
    Act                   &   dim             &   method   & $e_0$              & $e_1$              & $e_2$  & time(s)  & NoP   \\ 
\hline
\multirow{9}{*}{ReQU} & \multirow{3}{*}{d=2}  & Mix & \textbf{0.0003} & \textbf{0.0013} & \textbf{0.0006}  & 3.88   &5410\\  
\cline{3-8}
                      &                       & DGM   & 0.0024          & 0.0029          & 0.0083           & 5.17   &5589\\ 
\cline{3-8}
                      &                       & DRM  & 0.038           & 0.0308          & 0.07993          & 1.61   &5589\\ 
\cline{2-8}
                      & \multirow{3}{*}{d=5}  & Mix & \textbf{0.0091} & \textbf{0.0255} & \textbf{0.0054}  & 5.89   &32650\\ 
\cline{3-8}
                      &                       & DGM   & 0.0168          & 0.0376          & 0.0226           & 8.31   &33402\\ 
\cline{3-8}
                      &                       & DRM  & 0.0214          & 0.0301          & 0.1607           & 1.53   &33402\\ 
\cline{2-8}
                      & \multirow{3}{*}{d=10} & Mix & 0.0565          & 0.0845          & \textbf{0.0082} & 9.61   &129050 \\ 
\cline{3-8}
                      &                       & DGM   & \textbf{0.0463} & \textbf{0.0733} & 0.0326           & 13.31   &130063\\ 
\cline{3-8}
                      &                       & DRM  & 0.0737          & 0.0917          & 0.2608           & 1.72   &130063\\ 
\hline
\multirow{9}{*}{ReCU} & \multirow{3}{*}{d=2}  & Mix & \textbf{0.0001} & \textbf{0.0003} & \textbf{0.0003}  & 3.50   &5410\\ 
\cline{3-8}
                      &                       & DGM   & \textbf{0.0001} & \textbf{0.0003} & 0.0005           & 5.07   &5589\\ 
\cline{3-8}
                      &                       & DRM  & 0.0104          & 0.0077          & 0.0214          & 1.65   &5589 \\ 
\cline{2-8}
                      & \multirow{3}{*}{d=5}  & Mix & 0.0071          & 0.0209          & 0.0021           & 5.94   &32650\\ 
\cline{3-8}
                      &                       & DGM   & \textbf{0.0039} & \textbf{0.0089} & \textbf{0.0019}  & 8.47   &33402\\ 
\cline{3-8}
                      &                       & DRM  & 0.0255          & 0.0245          & 0.083            & 1.65   &33402\\ 
\cline{2-8}
                      & \multirow{3}{*}{d=10} & Mix & 0.04991         & 0.0626          & \textbf{0.0057}  & 9.80   &129050\\ 
\cline{3-8}
                      &                       & DGM   & \textbf{0.0234} & \textbf{0.0403} & 0.01             & 13.51   &130063\\ 
\cline{3-8}
                      &                       & DRM  & 0.066           & 0.0838          & 0.2387          & 1.86   & 130063\\
\hline
\end{tabular}
\caption{Errors of three neural networks methods with different activation functions (Act) on the Dirichlet problem. For dimesion (dim) $d=2$, $5$ and $10$, the $w$ is chosen to be $10$, $25$ and $50$ respectively. The averaged time per 100 iterations and number of parameters (NoP) are listed in the last two columns as well. }
\label{table:poisson_D}
\end{table}

We also test the elliptic PDE with Neumann boundary condition:
\begin{equation*}\begin{aligned}
    -\Delta u(x) + u(x) = (\pi^2+1)\sum_{i=1}^d\cos(\pi x_i) \label{eq:elliptic_N} \text{ in }\Omega,\quad \frac{\partial u(x)}{\partial\mathbf{n}} = 0 \text{ on } \partial\Omega.\notag
\end{aligned}\end{equation*}
The true solution is $u^*(x)=\sum_{i=1}^d\cos(\pi x_i)$. 
% The expected loss function for the mixed residual method is (\ref{eq:expected_loss_N}) with $\lambda_1=100$ and $\lambda_2=10$. The loss of DRM is:
% \begin{equation*}
%     \Vert\phi-(\pi^2+1)\sum_{i=1}^d\cos(\pi x_i)\Vert_{L^2(\Omega)}^2+\Vert\nabla\phi\Vert_{L^2(\Omega)}^2.
% \end{equation*}
% The loss of DGM is:
% \begin{equation*}
%     \Vert-\Delta\phi+\phi-(\pi^2+1)\sum_{i=1}^d\cos(\pi x_i)\Vert_{L^2(\Omega)}^2+\lambda_1\Vert\frac{\partial \phi}{\partial\mathbf{n}}-0\Vert_{L^2(\partial\Omega)}^2,
% \end{equation*}
%where $\lambda_1=100$. 
The relative errors are shown in Table \ref{table:poisson_N}. We also observe that the DRM has the largest relative errors in all the cases, and the mixed residual method achieve similar accuracy with the benchmark while it cost less time.

\begin{table}
\centering
\begin{tabular}{|p{0.08\textwidth}|p{0.07\textwidth}|p{0.1\textwidth}|p{0.1\textwidth}|p{0.1\textwidth}|p{0.1\textwidth}|p{0.09\textwidth}|p{0.08\textwidth}|} 
    \hline
        Act                   &   dim             &   method   & $e_0$              & $e_1$              & $e_2$  & time(s)  & NoP   \\ 
    \hline
\multirow{9}{*}{ReQU} & \multirow{3}{*}{d=2}  & Mix & \textbf{0.0006} & 0.002           & \textbf{0.0025} & 3.58   & 5410  \\ 
\cline{3-8}
                      &                       & DGM   & 0.0017          & \textbf{0.0017} & 0.0051    & 4.24   & 5589       \\ 
\cline{3-8}
                      &                       & DRM  & 0.077           & 0.0136          & 0.0847    & 1.99   & 5589       \\ 
\cline{2-8}
                      & \multirow{3}{*}{d=5}  & Mix & \textbf{0.0027} & 0.0096          & 0.0086    & 5.00   & 32650       \\ 
\cline{3-8}
                      &                       & DGM   & 0.0077          & \textbf{0.0027} & \textbf{0.0076} & 7.44   & 33402 \\ 
\cline{3-8}
                      &                       & DRM  & 0.0377          & 0.0484          & 0.12      & 1.85   & 33402       \\ 
\cline{2-8}
                      & \multirow{3}{*}{d=10} & Mix & \textbf{0.0114} & 0.0258          & \textbf{0.0187}  & 8.55   & 129050\\ 
\cline{3-8}
                      &                       & DGM   & 0.0123          & \textbf{0.016}  & 0.0253   & 12.76   & 130063        \\ 
\cline{3-8}
                      &                       & DRM  & 0.0446          & 0.0475          & 0.1439   & 2.10   & 130063        \\ 
\hline
\multirow{9}{*}{ReCU} & \multirow{3}{*}{d=2}  & Mix & \textbf{0.0000}      & \textbf{0.0002} & \textbf{0.0003}  & 3.86   & 5410\\ 
\cline{3-8}
                      &                       & DGM   & 0.0003          & \textbf{0.0002} & 0.0005   & 4.76   & 5589        \\ 
\cline{3-8}
                      &                       & DRM  & 0.0162          & 0.0019          & 0.0292    & 1.85   & 5589       \\ 
\cline{2-8}
                      & \multirow{3}{*}{d=5}  & Mix & 0.0006          & 0.0013          & \textbf{0.0009}  & 5.67   & 32650\\ 
\cline{3-8}
                      &                       & DGM   & \textbf{0.0005} & \textbf{0.0004} & \textbf{0.0009}  & 8.21   & 33402\\ 
\cline{3-8}
                      &                       & DRM  & 0.0176          & 0.0172          & 0.0324   & 1.68   & 33402        \\ 
\cline{2-8}
                      & \multirow{3}{*}{d=10} & Mix & \textbf{0.002}  & 0.004           & 0.0038   & 8.84   & 129050       \\ 
\cline{3-8}
                      &                       & DGM   & 0.0022          & \textbf{0.0016} & \textbf{0.0013}  & 13.16   & 130063\\ 
\cline{3-8}
                      &                       & DRM  & 0.0852          & 0.0252          & 0.0479   & 1.90   & 130063        \\
\hline
\end{tabular}
\caption{Errors of three neural networks methods with different activation functions (Act) on the Nuemman problem. For dimension (dim) $d=2$, $5$ and $10$, the $w$ is chosen to be $10$, $25$ and $50$ respectively. The averaged time per 100 iterations and number of parameters (NoP) are listed in the last two columns as well. }
\label{table:poisson_N}
\end{table}

\begin{figure}
     \centering
     \begin{subfigure}[b]{0.27\textwidth}
         \centering
         \includegraphics[width=\textwidth]{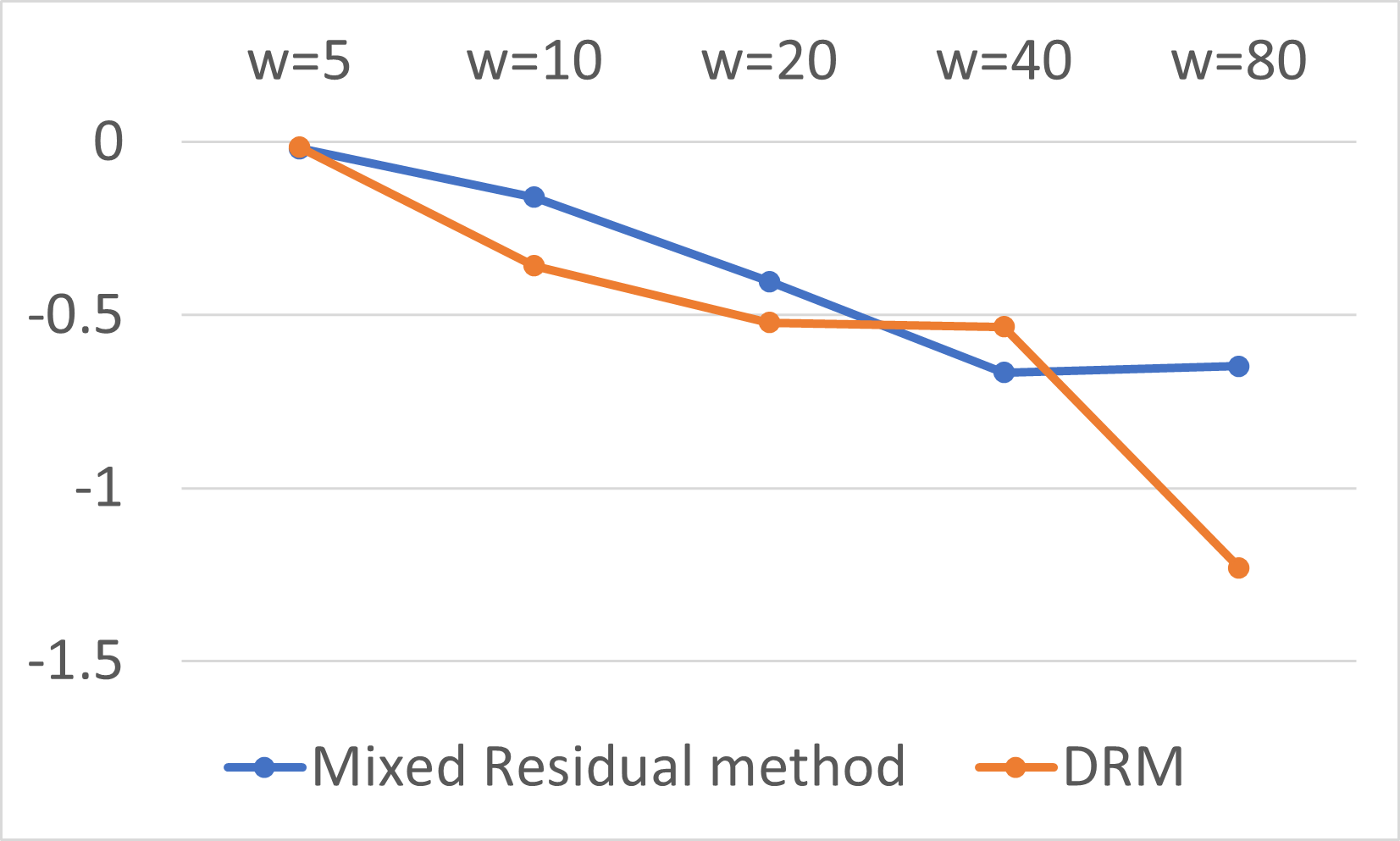}
         \caption{$\log(e_0)$}
     \end{subfigure}
     \hfill
     \begin{subfigure}[b]{0.27\textwidth}
         \centering
         \includegraphics[width=\textwidth]{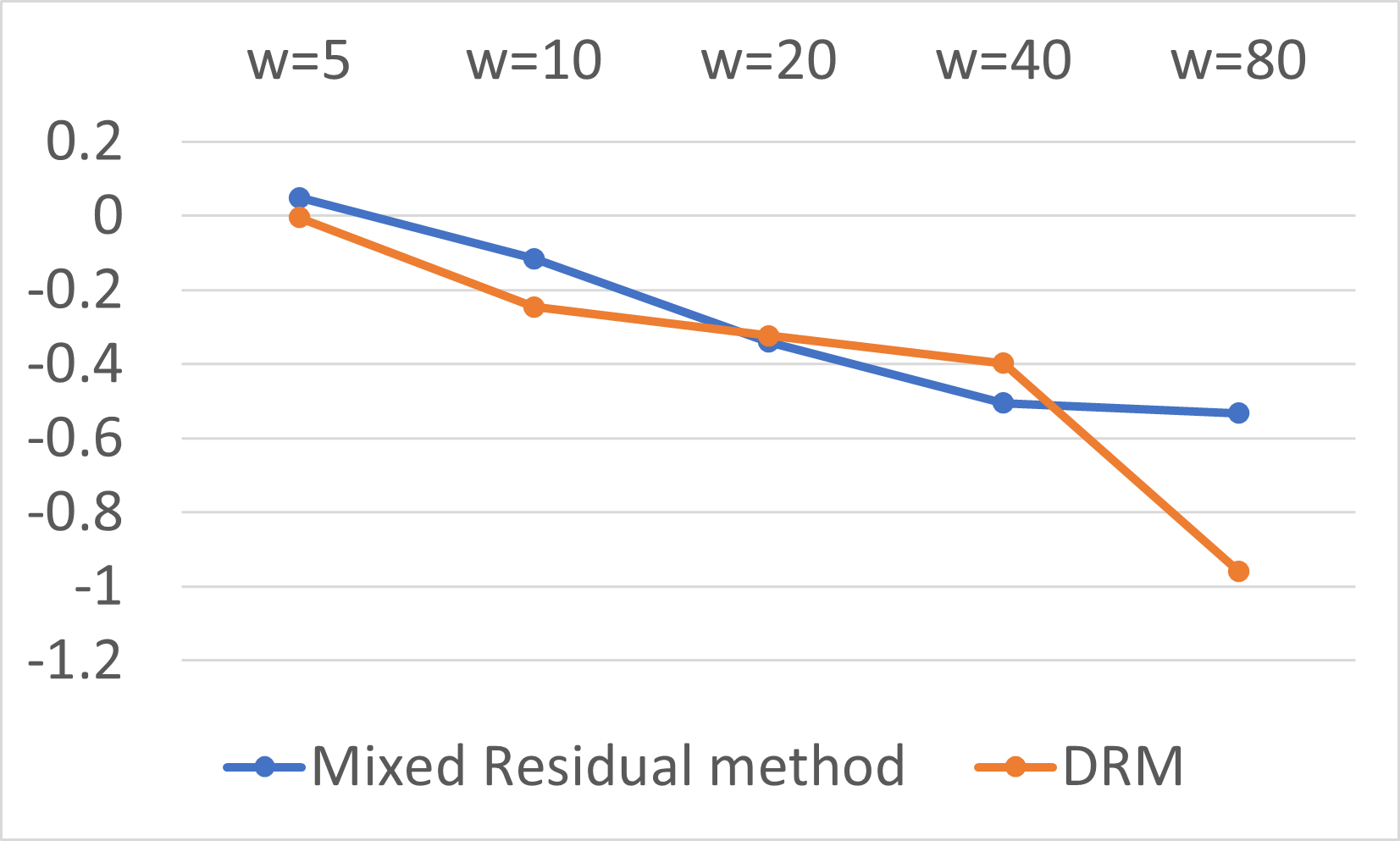}
         \caption{$\log(e_1)$}
     \end{subfigure}
     \hfill
     \begin{subfigure}[b]{0.27\textwidth}
         \centering
         \includegraphics[width=\textwidth]{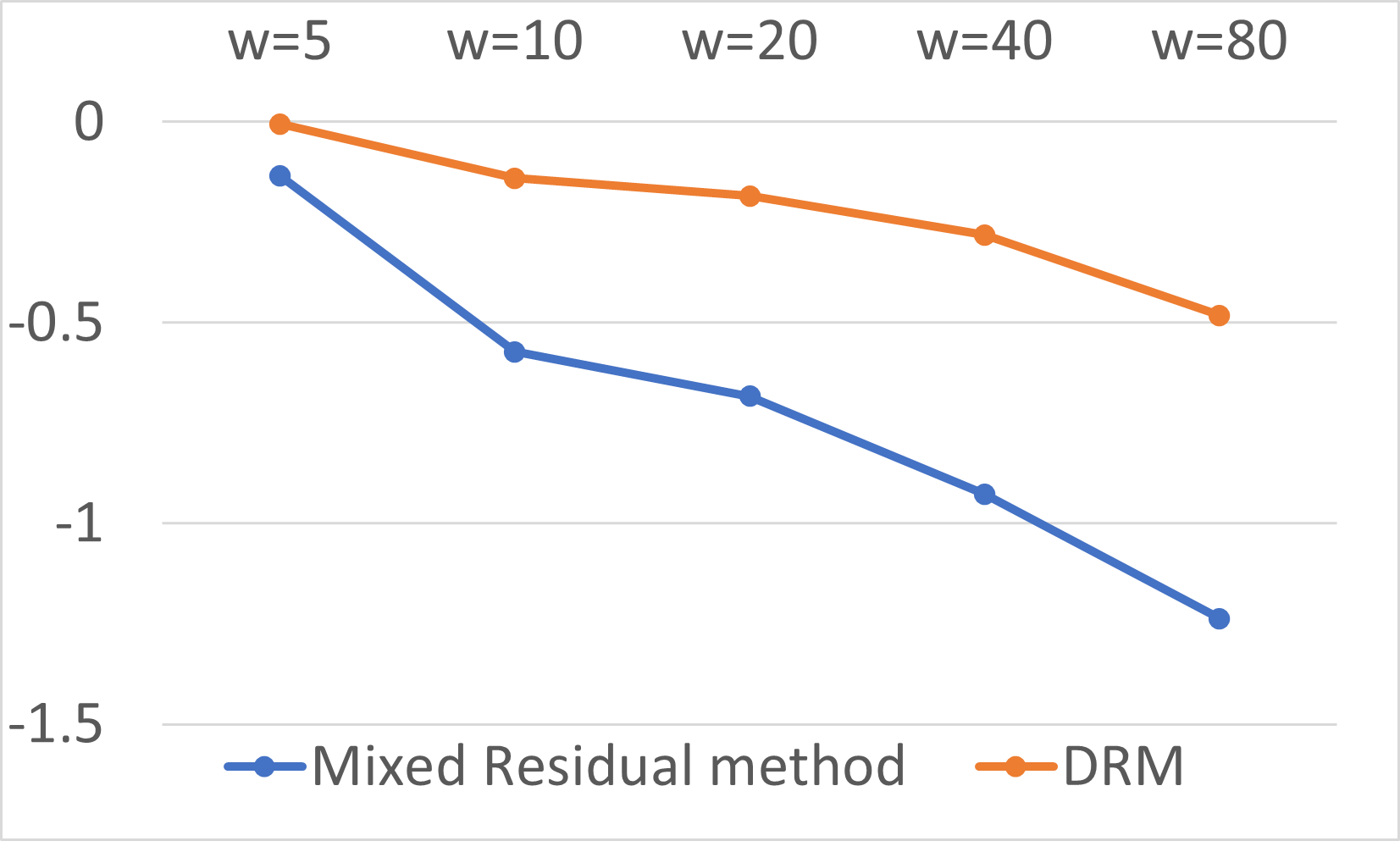}
         \caption{$\log(e_2)$}
     \end{subfigure}\\
     \begin{subfigure}[b]{0.27\textwidth}
         \centering
         \includegraphics[width=\textwidth]{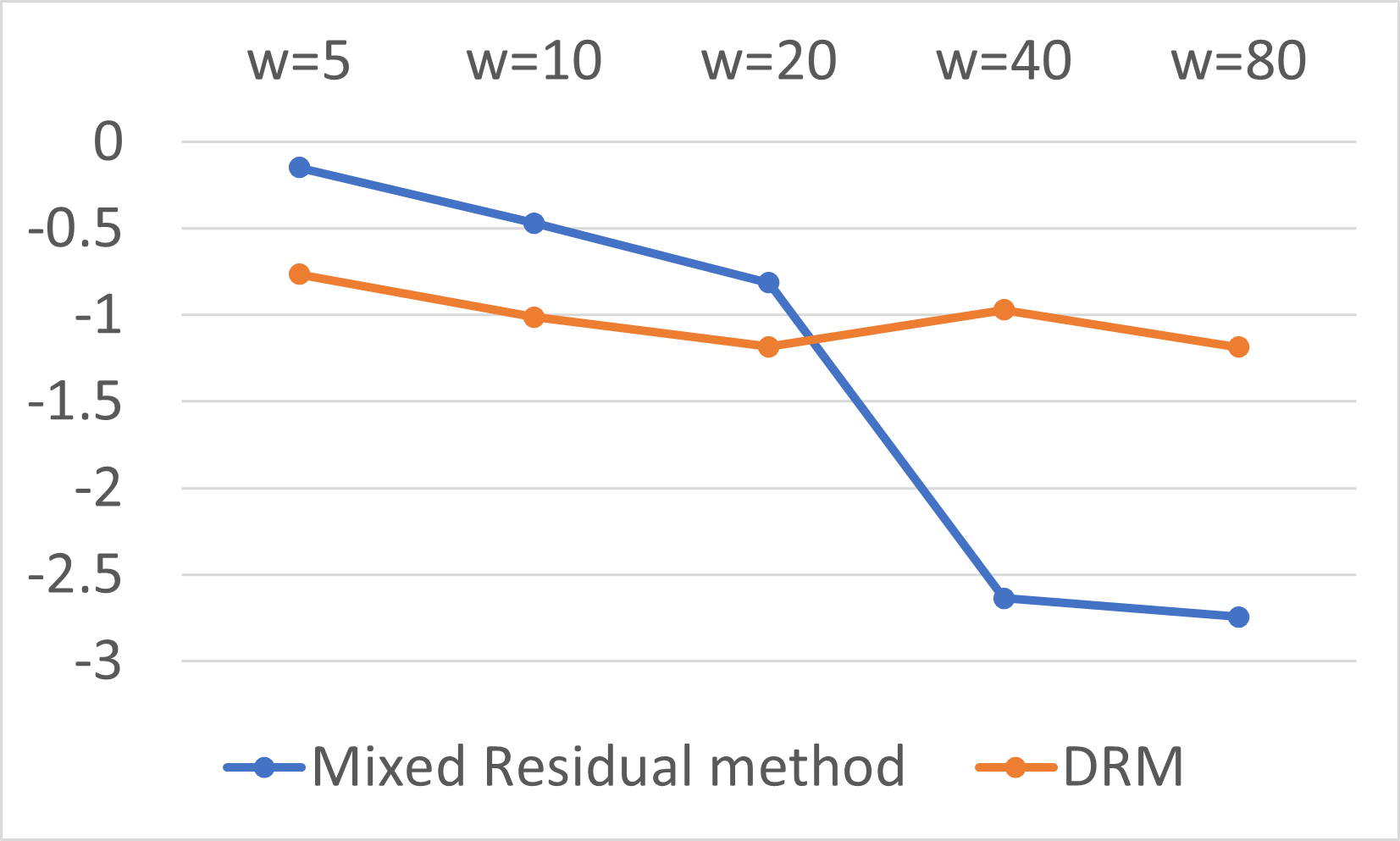}
         \caption{$\log(e_0)$}
     \end{subfigure}
     \hfill
     \begin{subfigure}[b]{0.27\textwidth}
         \centering
         \includegraphics[width=\textwidth]{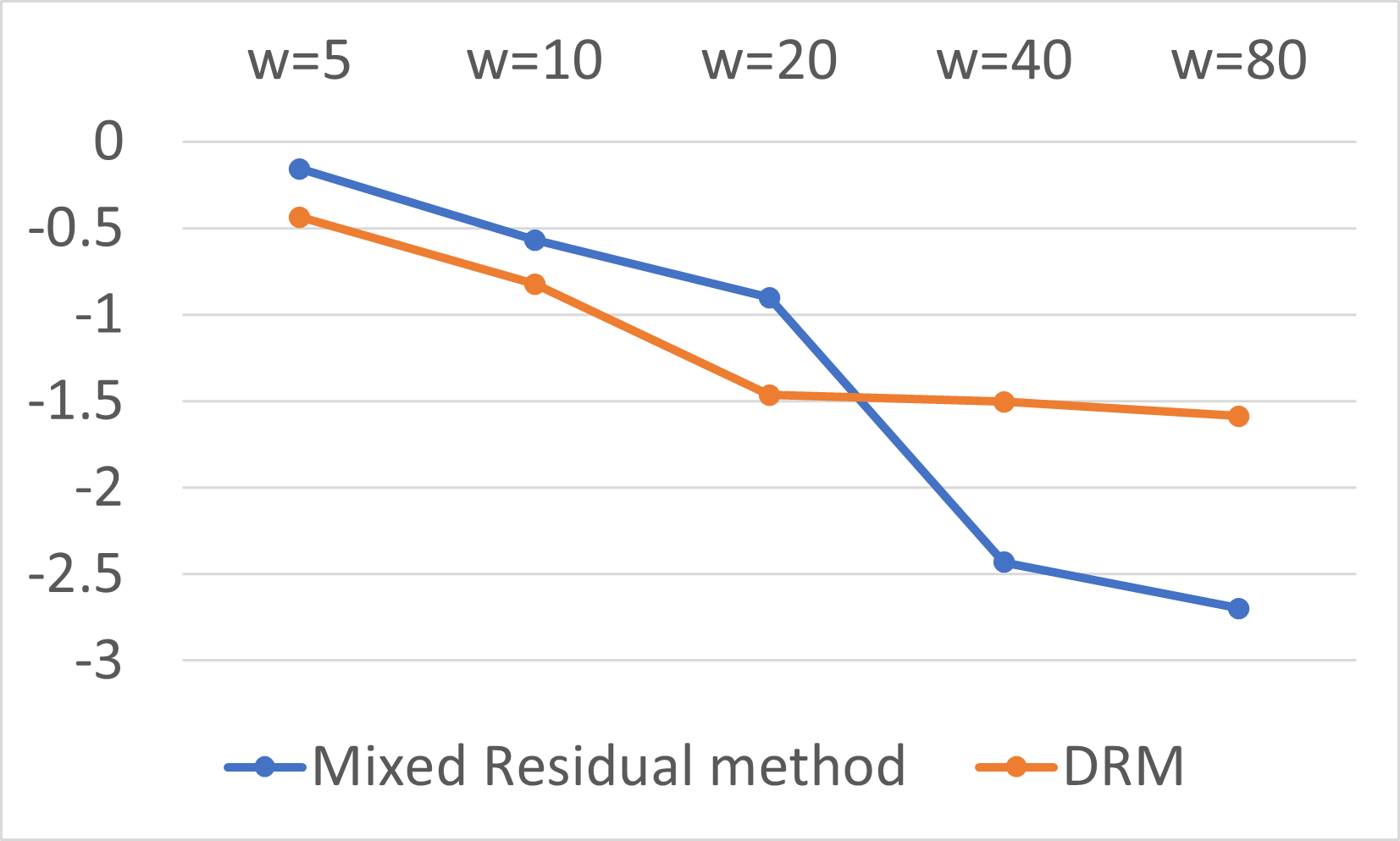}
         \caption{$\log(e_1)$}
     \end{subfigure}
     \hfill
     \begin{subfigure}[b]{0.27\textwidth}
         \centering
         \includegraphics[width=\textwidth]{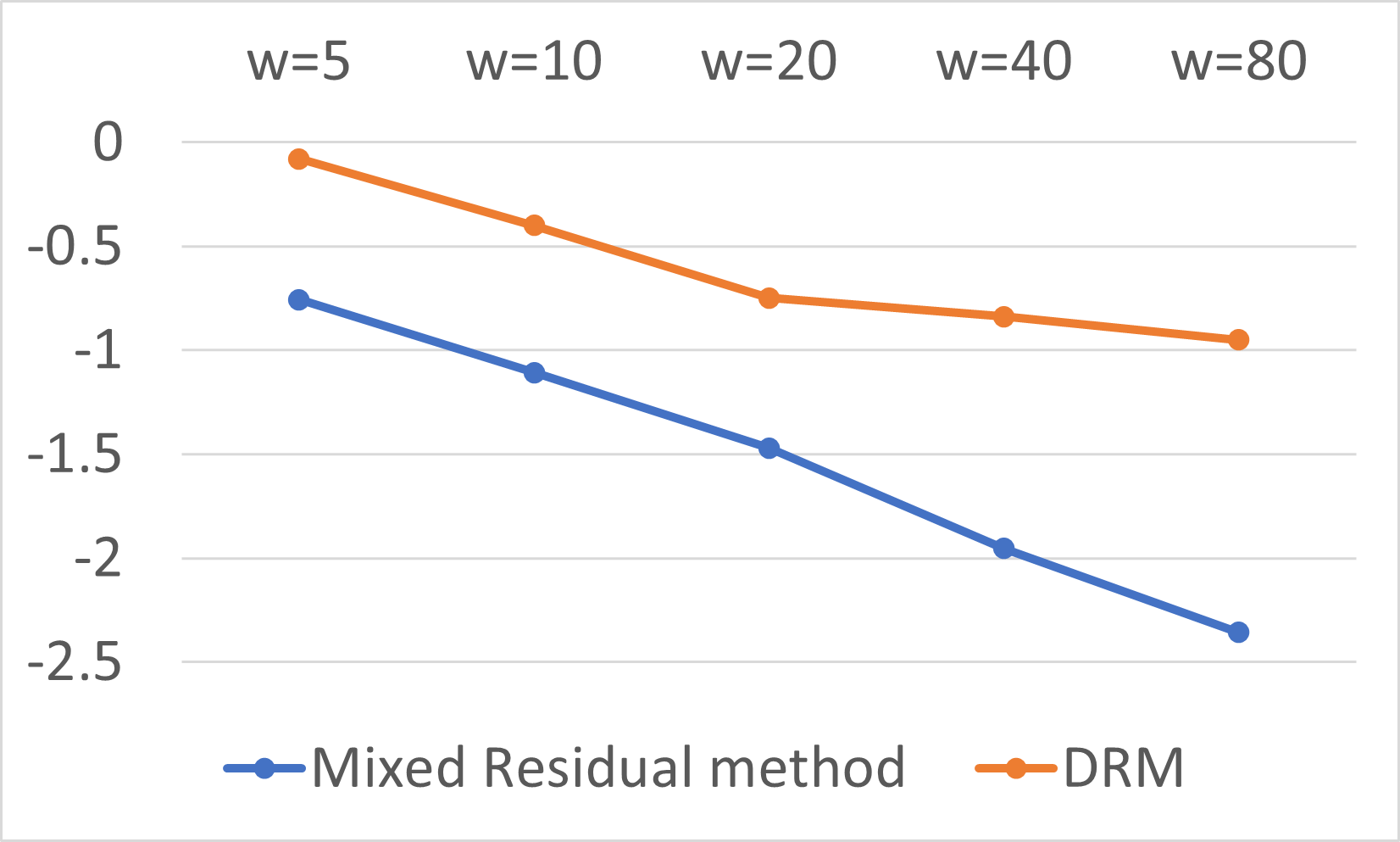}
         \caption{$\log(e_2)$}
     \end{subfigure}
        \caption{The relative errors of mixed residual method and DRM with different network size on the Dirichlet problem ((a)-(c)) and the Neumann problem ((d)-(f)). The depth of the networks are fixed to be $2$.}
        \label{fig:errors}
\end{figure}

We also want to compare the relative error of the mixed residual method and DRM when we increase the width of networks. We fix the depth to be 2, and set $w=5$, $10$, $20$, $40$ and $80$ respectively. Then, we evaluate the relative error on Dirichlet and Neumann problems with $d=10$. The logarithm of $e_0$, $e_1$ and $e_2$ after 500,000 iterations are shown in Figure \ref{fig:errors}. Generally, the $e_0$ and $e_1$ errors of DRM are competitive with the mixed residual method, while the $e_2$ error of the mixed residual method always decays much faster than DRM. The behavior of the two methods on the $e_2$ meets our expectation and verifies our analysis. 

\bibliographystyle{siamplain}
\bibliography{references}

\begin{thebibliography}{10}

\bibitem{barron1993universal}
{\sc A.~R. Barron}, {\em Universal approximation bounds for superpositions of a
  sigmoidal function}, IEEE Transactions on Information theory, 39 (1993),
  pp.~930--945.

\bibitem{cai2020deep}
{\sc Z.~Cai, J.~Chen, M.~Liu, and X.~Liu}, {\em Deep least-squares methods: An
  unsupervised learning-based numerical method for solving elliptic pdes},
  Journal of Computational Physics, 420 (2020), p.~109707.

\bibitem{cai1994first}
{\sc Z.~Cai, R.~Lazarov, T.~A. Manteuffel, and S.~F. McCormick}, {\em
  First-order system least squares for second-order partial differential
  equations: Part i}, SIAM Journal on Numerical Analysis, 31 (1994),
  pp.~1785--1799.

\bibitem{chen2018neural}
{\sc R.~T. Chen, Y.~Rubanova, J.~Bettencourt, and D.~Duvenaud}, {\em Neural
  ordinary differential equations}, arXiv preprint arXiv:1806.07366,  (2018).

\bibitem{gingold1977smoothed}
{\sc R.~A. Gingold and J.~J. Monaghan}, {\em Smoothed particle hydrodynamics:
  theory and application to non-spherical stars}, Monthly notices of the royal
  astronomical society, 181 (1977), pp.~375--389.

\bibitem{han2018solving}
{\sc J.~Han, A.~Jentzen, and E.~Weinan}, {\em Solving high-dimensional partial
  differential equations using deep learning}, Proceedings of the National
  Academy of Sciences, 115 (2018), pp.~8505--8510.

\bibitem{he2016deep}
{\sc K.~He, X.~Zhang, S.~Ren, and J.~Sun}, {\em Deep residual learning for
  image recognition}, in Proceedings of the IEEE conference on computer vision
  and pattern recognition, 2016, pp.~770--778.

\bibitem{ledoux1991probability}
{\sc M.~Ledoux and M.~Talagrand}, {\em Probability in Banach Spaces:
  isoperimetry and processes}, vol.~23, Springer Science \& Business Media,
  1991.

\bibitem{li2021generalization}
{\sc L.~Li, X.-C. Tai, and J.~Yang}, {\em Generalization error analysis of
  neural networks with gradient based regularization}, arXiv preprint
  arXiv:2107.02797,  (2021).

\bibitem{li2020fourier}
{\sc Z.~Li, N.~Kovachki, K.~Azizzadenesheli, B.~Liu, K.~Bhattacharya,
  A.~Stuart, and A.~Anandkumar}, {\em Fourier neural operator for parametric
  partial differential equations}, arXiv preprint arXiv:2010.08895,  (2020).

\bibitem{long2019pde}
{\sc Z.~Long, Y.~Lu, and B.~Dong}, {\em Pde-net 2.0: Learning pdes from data
  with a numeric-symbolic hybrid deep network}, Journal of Computational
  Physics, 399 (2019), p.~108925.

\bibitem{long2018pde}
{\sc Z.~Long, Y.~Lu, X.~Ma, and B.~Dong}, {\em Pde-net: Learning pdes from
  data}, in International Conference on Machine Learning, PMLR, 2018,
  pp.~3208--3216.

\bibitem{lu2021priori}
{\sc J.~Lu, Y.~Lu, and M.~Wang}, {\em A priori generalization analysis of the
  deep ritz method for solving high dimensional elliptic equations}, arXiv
  preprint arXiv:2101.01708,  (2021).

\bibitem{lu2019deeponet}
{\sc L.~Lu, P.~Jin, and G.~E. Karniadakis}, {\em Deeponet: Learning nonlinear
  operators for identifying differential equations based on the universal
  approximation theorem of operators}, arXiv preprint arXiv:1910.03193,
  (2019).

\bibitem{lu2021deepxde}
{\sc L.~Lu, X.~Meng, Z.~Mao, and G.~E. Karniadakis}, {\em Deepxde: A deep
  learning library for solving differential equations}, SIAM Review, 63 (2021),
  pp.~208--228.

\bibitem{lyu2020enforcing}
{\sc L.~Lyu, K.~Wu, R.~Du, and J.~Chen}, {\em Enforcing exact boundary and
  initial conditions in the deep mixed residual method}, arXiv preprint
  arXiv:2008.01491,  (2020).

\bibitem{lyu2022mim}
{\sc L.~Lyu, Z.~Zhang, M.~Chen, and J.~Chen}, {\em Mim: A deep mixed residual
  method for solving high-order partial differential equations}, Journal of
  Computational Physics,  (2022), p.~110930.

\bibitem{ma2019barron}
{\sc C.~Ma, L.~Wu, et~al.}, {\em Barron spaces and the compositional function
  spaces for neural network models}, arXiv preprint arXiv:1906.08039,  (2019).

\bibitem{massart2000some}
{\sc P.~Massart}, {\em Some applications of concentration inequalities to
  statistics}, in Annales de la Facult{\'e} des sciences de Toulouse:
  Math{\'e}matiques, vol.~9, 2000, pp.~245--303.

\bibitem{mishra2020estimates}
{\sc S.~Mishra and R.~Molinaro}, {\em Estimates on the generalization error of
  physics informed neural networks (pinns) for approximating pdes}, SAM
  Research Report, 2020 (2020).

\bibitem{mishra2020estimates2}
{\sc S.~Mishra and R.~Molinaro}, {\em Estimates on the generalization error of
  physics informed neural networks (pinns) for approximating pdes ii: A class
  of inverse problems}, SAM Research Report, 2020 (2020).

\bibitem{muller2021error}
{\sc J.~M{\"u}ller and M.~Zeinhofer}, {\em Error estimates for the variational
  training of neural networks with boundary penalty}, arXiv preprint
  arXiv:2103.01007,  (2021).

\bibitem{raissi2019physics}
{\sc M.~Raissi, P.~Perdikaris, and G.~E. Karniadakis}, {\em Physics-informed
  neural networks: A deep learning framework for solving forward and inverse
  problems involving nonlinear partial differential equations}, Journal of
  Computational Physics, 378 (2019), pp.~686--707.

\bibitem{siegel2020approximation}
{\sc J.~W. Siegel and J.~Xu}, {\em Approximation rates for neural networks with
  general activation functions}, Neural Networks, 128 (2020), pp.~313--321.

\bibitem{sirignano2018dgm}
{\sc J.~Sirignano and K.~Spiliopoulos}, {\em Dgm: A deep learning algorithm for
  solving partial differential equations}, Journal of computational physics,
  375 (2018), pp.~1339--1364.

\bibitem{weinan2018deep}
{\sc E.~Weinan and B.~Yu}, {\em The deep ritz method: a deep learning-based
  numerical algorithm for solving variational problems}, Communications in
  Mathematics and Statistics, 6 (2018), pp.~1--12.

\bibitem{yang2021local}
{\sc J.~Yang and Q.~Zhu}, {\em A local deep learning method for solving high
  order partial differential equations}, NUMERICAL MATHEMATICS-THEORY METHODS
  AND APPLICATIONS,  (2021).

\bibitem{zang2020weak}
{\sc Y.~Zang, G.~Bao, X.~Ye, and H.~Zhou}, {\em Weak adversarial networks for
  high-dimensional partial differential equations}, Journal of Computational
  Physics, 411 (2020), p.~109409.

\end{thebibliography}

\appendix
\section{Proof of lemma (\ref{theorem:H1sp})}\label{app_approx_er}
The proof starts with the lemma 4.5 in \cite{lu2021priori}:
\begin{lemma}\cite{lu2021priori}
    Let $ \Omega = [-1,1],$ $g\in C^2(\Omega)$ with $\|g^{(s)}\|_{L^\infty(\Omega)}\leq B$ for $s = 0,1,2$. Assume that $g'(0)=0$. Let $\{z_i\}_{i=0}^{2m}$ be a uniform mesh on $\Omega$, where $h = \frac{1}{m}$ and \begin{equation*}
        z_i = -1 + ih, i =0,\cdots, 2m.
        \label{eq::partition}
    \end{equation*} 
    Then there exists a two-layer ReLU neural network $g_m(z)$ of the form 
    \begin{equation}
        g_m(z)  = c + \sum_{i=1}^{2m} a_i \text{\rm ReLU}(\epsilon_i z - b_i), z\in\Omega,
        \label{eq::origin_func}
    \end{equation}
    with $c=g(0), b_i\in[-1,1], $ and $\epsilon_i \in\{\pm 1\},$ $i=1\cdots,2m$ such that 
    \begin{equation*}
        \|g - g_m\|_{W^{1,\infty}(\Omega)}\leq \frac{2B}{m}.
    \end{equation*}
    \end{lemma}
In addition, an example of the coefficients $\{a_i\}_{i=1}^{2m}$ are given as 
\begin{equation}
    a_i = \left\{\begin{aligned}
        &\frac{g(z_{m+1})-g(z_m)}{h}, & i = m + 1,\\
        &\frac{g(z_{m-1})-g(z_m)}{h}, & i = m,\\
        &\frac{g(z_i)-2g(z_{i-1})+g(z_{i-2})}{h}, &i>m+1,\\
        &\frac{g(z_{i-1})-2g(z_i)+g(z_{i+1})}{h}, &i<m.
    \end{aligned}\right.
    \label{coe::a_origin}
\end{equation}
and $\sum_{i=1}^{2m}|a_i|\leq 4B$, $|c|\leq B$.
This constructed function $g_m(z)$ is equalvant to a piecewise linear interpolation function on $\{z_i\}_{i=1}^{2m}$ and 
\begin{equation*}
    g_m(z) = g(0) + \sum_{i=1}^{2m}a_i\text{ReLU}(\epsilon_i(z-z_i)),  
\end{equation*} 
with $\epsilon_i = 1$ for $i = 1,\cdots,m$ and $\epsilon_i = -1$ for $i = m+1,\cdots,2m$. Notice that 
\begin{equation*}
    \begin{aligned}
        4\delta\relu(z) = \requ(z+\delta)-\requ(z-\delta) + e(z,\delta), 
    \end{aligned}
\end{equation*}
and 
\begin{equation*}
     e(z;\delta) = \left\{\begin{aligned}
            % &0, &z \leq -\delta,\\
            &-(z + \delta)^2, &-\delta< z\leq 0,\\
            &-(z - \delta)^2, & 0< z \leq \delta,\\
            &0, &\text{otherwise}.
        \end{aligned}\right.
        \label{eq::e_delta}
\end{equation*}
For fixed $\delta = h$, it follows that $|e(z)|\leq h^2$ and $|e'(z)|\leq 2h$.
Then we obtain a ReQU activated neural network $\hat{g}_m(z)$ by approximating the ReLU active function in (\ref{eq::origin_func}):
\begin{equation}
    \hat{g}_m(z) = g(0) + \sum_{i=1}^{2m}\frac{a_i}{4h}(\text{ReQU}(\epsilon_i(z-z_i)+h)-\text{ReQU}(\epsilon_i(z-z_i)-h)).
    \label{eq::relu2_origin}
\end{equation} 
Since $e(z)\neq 0$ only for $z\in [-h,h]$, it follows directly that 
    \begin{equation}
        \begin{aligned}
            \|g - \hat{g}_m\|_{W^{1,\infty}(\Omega)}&\leq  \|g - g_m\|_{W^{1,\infty}(\Omega)} + \|g_m - \hat{g}_m\|_{W^{1,\infty}(\Omega)}\\
            &\leq \frac{2B}{m} + \frac{3B}{m} = \frac{5B}{m}.
        \end{aligned}
        \label{eq::estimation_g}
    \end{equation}
Rearrange (\ref{eq::relu2_origin}) in the form of 
    \begin{equation*}
            \hat{g}_m(z) = g(0) + \sum_{i=0}^{2m+3} \hat{a}_i \requ(\epsilon_i z - b_i),
    \end{equation*}
    where 
    \begin{equation*}
        \hat{a}_i =\left\{\begin{aligned}
            &-\frac{a_{i+1}}{4h}, && i = 0,1,\\
            &\frac{a_{i-1} - a_{i+1}}{4h}, && 2\leq i \leq m-1,\\
            &\frac{a_{i-1}}{4h}, && m\leq i \leq m+3,\\
            &\frac{a_{i-1} - a_{i-3}}{4h}, && m+4\leq i \leq 2m+1,\\
            &-\frac{a_{i-3}}{4h}, && i = 2m+2,2m+3,
        \end{aligned}\right.
    \end{equation*}
    and
        \begin{equation*}
        b_i = \left\{\begin{aligned}
            &\epsilon_i z_i, && 0\leq i \leq m+1,\\
            &\epsilon_i z_{i-3}, && m+2\leq i \leq 2m+3.
        \end{aligned}\right.
    \end{equation*}
    Recall the formualtion (\ref{coe::a_origin}) of $a_i$. For $i\in [2,m-1]\cup[m+4,2m+1]$, we have 
    \begin{equation*}
        \begin{aligned}
            |\hat{a}_i|  &= \frac{|a_{s+1} - a_{s-1}|}{4h} (s=i\text{ or } i-2) \\
                & \leq \frac{2h}{4h} |g''(\xi_1) - g''(\xi_2)|\\
                & \leq |g^{(3)}(\xi_3)| \frac{|\xi_1 - \xi_2|}{2}\leq 2hB.
        \end{aligned}
    \end{equation*}
    Notice that the last step is satisfied if $g\in C^3(\Omega)$, otherwise $\sum_{i=0}^{2m+3} |\hat{a}_i|\sim O(m)$ and it is unacceptable. When $g\in C^3(\Omega)$, it concludes that 
        \begin{equation*}
        \begin{aligned}
            \sum_{i=0}^{2m+3} |\hat{a}_i|  \leq 4(m-1)hB + \frac{8}{4h}\cdot 2hB\leq 8B. 
        \end{aligned}
    \end{equation*}
    We have shown the approximation capability of the ReQU network, and now we further associate it with Barron space. 
    The second step can be completed by lemma 1 in \cite{barron1993universal} and theorem 4.1 in \cite{lu2021priori}:
    \begin{lemma}\cite{barron1993universal}
    Let $\mathcal{G}$ be a set in a Hilbert space, and $u$ lies in the closure of the convex hull of $\mathcal{G}$. Every element in $\mathcal{G}$ is bounded, i.e., $\forall g\in \mathcal{G}$, $\|g\|\leq B$. Then for every $m\in\mathbb{N}$, there exists $\{g_i\}_{i=1}^m\subset \mathcal{G}$ and $\{\lambda_i\}_{i=1}^m\subset[0,1]$ with $\sum_{i=1}^m\lambda_i = 1$ such that 
    \begin{equation}
        \|u -  \sum_{i=1}^m\lambda_ig_i\| \leq \frac{B}{\sqrt{m}}.
    \end{equation}
    \label{lem::convexhull}
    \end{lemma}
   \begin{lemma}\cite{lu2021priori}
    Let $u\in \mathcal{B}(\Omega)$. Then there exists $u_m$ which is a convex combination of $m$ functions in $\mathcal{F}_{\cos}(B)$ with $B = \|u\|_{\mathcal{B}(\Omega)}$ such that 
    \begin{equation*}
        \|u - \hat{u}(0) - u_m\|_{H^1(\Omega)}\leq \frac{\|u\|_{\mathcal{B}(\Omega)}}{\sqrt{m}}.
    \end{equation*}
\end{lemma}
Different from $u\in\mathcal{B}(\Omega)=\mathcal{B}^2(\Omega)$ in \cite{lu2021priori}, we let $u\in\mathcal{B}(\Omega)=\mathcal{B}^3(\Omega)$. Accordingly, $\mathcal{F}_{\cos}(B)$ should be changed to 
\begin{equation*}
    \mathcal{F}_{\cos}(B):=\{\frac{\gamma}{1+\pi^3|k|^3}\cos(\pi(k\cdot x + b))|k\in\mathbb{Z}^d \backslash \{\textbf{0}\}, |\gamma|\leq B,b\in\{0,1\}\}.
\end{equation*}
and the result still holds.      

Notice that every function in $\mathcal{F}_{\cos} (B)$ is the composition of the one dimensional function $g$ defined on $[-1,1]$ by 
\begin{equation*}
    g(z) = \frac{\gamma}{1+\pi^3|k|^3}\cos(\pi(|k|z+b)),
\end{equation*}
with $k\in\mathbb{Z}^d\backslash\{\textbf{0}\},|\gamma|\leq B$ and $b\in\{0,1\}$, and a linear function $z = \omega\cdot x$ with $w = \frac{k}{|k|}$. It is clear that $g \in C^3([-1,1])$ satisfies that 
\begin{equation*}
    \|g^{(s)}\|_{L^\infty([-1,1])}\leq |\gamma|\leq B \text{ for } s=0,1,2,3,
\end{equation*}
and $g'(0) =0$. 
Formulation (\ref{eq::estimation_g}) yields that 
\begin{equation*}
    \|g(\omega\cdot x)-\hat{g}(\omega\cdot x)\|_{H^1(\Omega)}\leq \|g-\hat{g}\|_{W^{1,\infty}(\Omega)}\leq \frac{5B}{m}.
\end{equation*}
Then $u$ lies in the closure of convex hull of 
\begin{equation*}
    \mathcal{F}_{\text{Q}}(B):= \{c + a\text{ReQU}(\omega^Tx+b)\vert |c|\leq 2B, |a|\leq 8B, \|w\|_2\leq 1, |b|\leq 1\}.
\end{equation*}
For any $u_f\in\mathcal{F}_{\text{Q}}(B)$, 
\begin{equation}
    \|u_f\|_{H^1(\Omega)}\leq (8d+16\sqrt{d}+10)B + 16(\sqrt{d}+1)B.
    \label{bound::F_Q}
\end{equation}
Thanks to lemma \ref{lem::convexhull} and the bound (\ref{bound::F_Q}), there exists a $u_m\in V^m_u$, which is a convex combination of $m$ functions in $\mathcal{F}_{\text{Q}}(\normu)$ such that
\begin{equation*}
    \|u - u_m\|_{H^1(\Omega)}\leq \frac{8d+32\sqrt{d}+26}{\sqrt{m}}\normu.
\end{equation*}
It completes the proof. There is a remark that for the purpose of making the coefficient bounded, we let $u\in \mathcal{B}^3
(\Omega)$, which is a sufficient condition. There may be weaker conditions or stronger conclusions here.

\section{Proof of lemma \ref{lem::computing_rules}} \label{app_CR}
    The proof of (i) and (ii) is trivial by the definition. A rigorous proof of (iv) can be found in \cite{ledoux1991probability} and (v) is a direct corollary of (iv) by letting $\sigma(x)=x^2$. We only prove the remaining properties.

(iii)
First use the definition of the Rademacher complexity 
\begin{equation*}
    \begin{aligned}
        R_n(g) = \mathbb{E}_X\mathbb{E}_\varepsilon\sup_{g\in\{g\}}|\frac{1}{n}\sum_{i=1}^n \varepsilon_i g(X_i)|= \mathbb{E}_X\mathbb{E}_\varepsilon|\frac{1}{n}\sum_{i=1}^n \varepsilon_i g_i|,
    \end{aligned}
\end{equation*}
where $g_i = g(X_i)$, $i= 1,\cdots, n$. Next we expand part of this summation 
\begin{equation}
    \begin{aligned}
        (\sum_{i=1}^n \varepsilon_i g_i)^2= \sum_{i=1}^n \varepsilon_i^2 g_i^2 + \sum_{i\neq j} \varepsilon_i\varepsilon_j g_ig_j\leq b^2(\sum_{i=1}^n\varepsilon_i)^2 + \sum_{i\neq j} \varepsilon_i\varepsilon_j (g_ig_j-b^2).\\
    \end{aligned}
\end{equation} 
Taking the expectation of $\varepsilon$ yields 
\begin{equation*}
    \begin{aligned}
        \mathbb{E}(\sum_{i=1}^n \varepsilon_i g_i)^2 & \leq b^2\mathbb{E}(\sum_{i=1}^n\varepsilon_i)^2 + \mathbb{E}[\sum_{i\neq j} \varepsilon_i\varepsilon_j (g_ig_j-b^2)]\\
        &=b^2\mathbb{E}(\sum_{i=1}^n\varepsilon_i)^2 + \sum_{i\neq j} \mathbb{E}[\varepsilon_i\varepsilon_j] (g_ig_j-b^2)\\
        &=b^2\mathbb{E}(\sum_{i=1}^n\varepsilon_i)^2 + \sum_{i\neq j} \mathbb{E}[\varepsilon_i]\mathbb{E}[\varepsilon_j] (g_ig_j-b^2)\\
        &=b^2\mathbb{E}(\sum_{i=1}^n\varepsilon_i)^2.
    \end{aligned}
    \label{pf::Eb2}
\end{equation*}
Let $Y = \sum_{i=1}^n\varepsilon_i$, and it satisfies
\begin{equation*}
    \mathbb{E}[Y]  = \sum_{i=1}^n\mathbb{E}[\varepsilon_i] = 0, \ \text{Var}[Y] =  n\text{Var}[\varepsilon] = n.
\end{equation*}
With the above estimation, the proof is completed 
\begin{equation*}
    \begin{aligned}
        R_n(g)\leq \frac{1}{n} \sqrt{\mathbb{E}(\sum_{i=1}^n \varepsilon_i g_i)^2}\leq  \frac{b}{n} \sqrt{\mathbb{E}[Y^2]} = \frac{b}{n}\sqrt{\text{Var}[Y]} = \frac{b}{\sqrt{n}}.
    \end{aligned}
\end{equation*}

(vi) 
Denote $f_i = f(X_i)$, $g_i = g(X_i)$ for $i = 1,\cdots,n$.
\begin{equation*}
    \begin{aligned}
        R_n(\mathcal{F}\mathcal{G})&= \mathbb{E}_{X,\varepsilon}\sup_{f\in\mathcal{F},g\in\mathcal{G}}|\frac{1}{n}\sum_{i=1}^n\varepsilon_if_ig_i|\\
        &=\mathbb{E}_{X,\varepsilon}\sup_{f\in\mathcal{F},g\in\mathcal{G}}|\frac{1}{n}\sum_{i=1}^n\varepsilon_i(\frac{1}{2}(f_i + g_i)^2 - \frac{1}{2}f_i^2-\frac{1}{2}g_i^2)|\\
        &\leq \frac{1}{2}R_n((\mathcal{F}+\mathcal{G})^2) + \frac{1}{2}R_n(\mathcal{F}^2) + \frac{1}{2}R_n(\mathcal{G}^2).
    \end{aligned}
\end{equation*}
Using the result of (v) yields 
\begin{equation*}
    \begin{aligned}
        R_n(\mathcal{F}\mathcal{G})&\leq 2\sup_{f\in\mathcal{F},g\in\mathcal{G}}\|f+g\|_\infty R_n(\mathcal{F}+\mathcal{G}) \\
        & + 2\sup_{f\in\mathcal{F}}\|f\|_\infty R_n(\mathcal{F}) + 2\sup_{g\in\mathcal{G}}\|g\|_\infty R_n(\mathcal{G})\\
        &\leq 4\sup_{f\in\mathcal{F}\cup \mathcal{G}}\|f\|_\infty R_n(\mathcal{F}+\mathcal{G}) +2\sup_{f\in\mathcal{F}\cup \mathcal{G}}\|f\|_\infty (R_n(\mathcal{F})+R_n(\mathcal{G}))\\
        &\leq 6\sup_{f\in\mathcal{F}\cup \mathcal{G}}\|f\|_\infty (R_n(\mathcal{F})+R_n(\mathcal{G})).
    \end{aligned}
\end{equation*}

\end{document}